\documentclass[a4paper, 11pt]{article}

\usepackage[utf8]{inputenc}
\usepackage[english]{babel}
\usepackage{amsthm, amsmath, amssymb}
\usepackage[T1]{fontenc}
\usepackage{fullpage}
\usepackage{lipsum}
\usepackage{authblk}
\usepackage{url} 
\usepackage{hyperref}
\usepackage{graphicx}
\usepackage{mathtools}
\usepackage{bbm}
\usepackage{tikz-cd,tikz-3dplot}
\usepackage{enumitem}
\usepackage{soul}
\usepackage{mathrsfs}
\usepackage{subcaption}
\usepackage{algorithm}
\usepackage{algpseudocode}

\usepackage[maxbibnames=99, style=numeric]{biblatex} 

\newtheorem{assumption}{Assumption}[section]
\newtheorem{definition}[assumption]{Definition}
\newtheorem{theorem}[assumption]{Theorem}
\newtheorem{proposition}[assumption]{Proposition}
\newtheorem{corollary}[assumption]{Corollary}

\theoremstyle{remark}

\theoremstyle{remark}
\newtheorem{remark}[assumption]{Remark}

\newcommand{\Ecal}{\mathcal{E}}

\newcommand{\Gcal}{\mathcal{G}}

\newcommand{\Ical}{\mathcal{I}}

\newcommand{\Scal}{\mathcal{S}}

\newcommand{\Wcal}{\mathcal{W}}

\newcommand{\Prob}{\mathscr{P}}
\newcommand{\Mart}{\mathscr{M}}
\newcommand{\Law}{{\mathscr{L}}}

\newcommand{\EE}{\mathbb{E}}

\newcommand{\NN}{\mathbb{N}}

\newcommand{\PP}{\mathbb{P}}
\newcommand{\QQ}{\mathbb{Q}}
\newcommand{\RR}{\mathbb{R}}

\newcommand{\RRd}{{\mathbb{R}^d}}

\newcommand{\mean}{{\textup{mean}}}

\newcommand{\supp}{{\textup{supp}}}
\newcommand{\co}{{\textup{co}}}

\newcommand{\MCov}{{\mathrm{MCov}}}
\newcommand{\interior}{{\mathrm{int}}}
\newcommand{\dom}{{\mathrm{dom}}}

\numberwithin{equation}{section}

\makeatletter
\renewcommand{\p@enumii}{}           

\makeatother

\makeatletter
\newenvironment{breakablealgorithm}
{%
\vspace{3mm}
\refstepcounter{algorithm}
\hrule height .8pt depth 0pt \kern 2pt
\renewcommand{\caption}[2][\relax]{
  {\raggedright\textbf{\ALG@name~\thealgorithm} ##2\par}
  \ifx\relax##1\relax
    \addcontentsline{loa}{algorithm}{\protect\numberline{\thealgorithm}##2}
  \else
    \addcontentsline{loa}{algorithm}{\protect\numberline{\thealgorithm}##1}
  \fi
  \kern 2pt\hrule\kern 2pt
}}
{
\kern 2pt\hrule
\vspace{3mm}
}
\makeatother

\title{$q_x$-Bass Martingales: Optimality, Existence and Convergence}
\author{Antonio Marini}

\affil{Department of Mathematics,\\ ETH Z\"{u}rich, Switzerland \\ \texttt{antonio.marini@math.ethz.ch}}
\begin{document}
\maketitle

\begin{abstract}
    We introduce the class of $q_x$-Bass martingales, extending the notion of $q$-Bass martingales. These are martingale couplings with prescribed marginals whose transition kernels are closest to a given family of reference kernels $(q_x)_{x \in \RRd}$ in the sense of a suitable weak optimal transport problem. We show that $q_x$-Bass martingales are optimizers of this problem and establish a strong duality result. Building on this duality, we introduce an algorithm for their computation and, under compactness assumptions on the prescribed marginals, prove its convergence, thereby establishing existence. We then specialize to one-dimensional marginals and show that the algorithm reduces to a fixed-point iteration.
\end{abstract}

\section{Introduction}

Over the past few years, martingale optimal transport has emerged as a powerful framework for constructing martingales with prescribed initial and terminal distributions while keeping their dynamics as close as possible to a given reference model. The existence of admissible martingales is guaranteed by Strassen’s theorem whenever the marginals are in convex order. Yet the theorem does not address the selection problem: it provides neither an explicit construction nor a criterion for distinguishing one martingale coupling among the admissible ones.

One prominent approach to this selection problem is based on Bass-type constructions, whose origin dates back to Richard Bass’ solution to the Skorokhod embedding problem \cite{Ba83}. Given a target distribution $\nu$, Bass’ approach starts by mapping a Gaussian random variable monotonically onto $\nu$ and then uses conditional expectations to construct a martingale with deterministic initial value and terminal law $\nu$.

This idea later inspired Backhoff-Veraguas--Beiglböck--Huesmann--Källblad \cite{BaBeHuKa20} to introduce the notion of stretched Brownian motion. This process is the optimizer of the martingale Benamou--Brenier formula, a martingale optimal transport problem which seeks a martingale with prescribed initial and terminal marginals whose dynamics are as close as possible to those of Brownian motion. Instead of considering a Brownian motion with deterministic initial value, they consider a Brownian motion with a suitable nontrivial initial law. By transporting the marginal of this Brownian motion at terminal time $T$ to the target distribution $\nu$ via the gradient of a convex function, and then taking conditional expectation, they construct a martingale whose initial law is a prescribed distribution $\mu$, dominated by $\nu$ in convex order. Given its close connection with Bass' construction, this martingale is often referred to in the literature as the Bass martingale from $\mu$ to $\nu$, and the initial law of the underlying Brownian motion is called the Bass distribution.

In the one-dimensional setting, Conze and Henry-Labordère \cite{CoHe21} introduced a numerical approach for determining the Bass distribution required to construct a Bass martingale from $\mu$ to $\nu$. Their method was developed for the calibration of the Bass local volatility model, which is based on this construction. In particular, their numerical method relies on a fixed-point iteration and, therefore, the Bass distribution is also often referred to as fixed-point distribution. The convergence of this method was studied by Acciaio--Marini--Pammer in \cite{AcMaPa23}, and a multidimensional extension of this algorithm was subsequently developed and studied by Hasenbichler--Joseph--Loeper--Ob{\l}oj--Pammer in \cite{HaJoLoObPa26}. 

In discrete time, Bass martingales were further generalized to an arbitrary reference measure $q$ by Tschiderer in \cite{Ts24}. Tschiderer showed that an analogous construction arises as the optimizer of a weak martingale optimal transport problem that selects a martingale coupling $(M_0,M_1)$ such that $M_0\sim\mu$, $M_1\sim\nu$, and the conditional laws $\Law(M_1| M_0=x)$, taken over all $x\in\supp(\mu)$, are as close as possible to the reference measure $q$ in a suitable sense. The properties of this class of discrete martingales were further studied by Acciaio and Marini in the one-dimensional setting; see \cite{AcMa26_semidiscrete, AcMa26_general}.

The present work further generalizes Tschiderer's framework. Given the initial and terminal marginals $\mu$ and $\nu$, and a reference kernel $(q_x)_{x\in\supp(\mu)}$, we introduce the following weak martingale optimal transport problem:
\begin{equation}
\label{eq:mainPb}
	WT_{\Ical}^{q_\cdot}(\mu,\nu) := \inf_{\pi\in\Mart(\mu,\nu)} \Ical(\pi), \qquad\qquad \Ical(\pi) := \int_\RRd \Wcal_2^2(\pi_x,q_x) \,\mu(dx),
\end{equation}
where $\Wcal_2$ denotes the $2$-Wasserstein distance and $\Mart(\mu,\nu)$ is the set of all martingale couplings with first marginal $\mu$ and second marginal $\nu$. Writing $\pi_x$ for the corresponding transition kernel, the martingale condition on the elements of $\Mart(\mu,\nu)$ is given by $\mean(\pi_x)=x$ for $\mu$-almost every $x$.

Under suitable integrability assumptions, the optimization problem \eqref{eq:mainPb} is equivalent to 
\begin{equation}
\label{eq:mainPb2}
WT_{\Scal}^{q_\cdot}(\mu,\nu) := \sup_{\pi\in\Mart(\mu,\nu)} \Scal(\pi), \qquad \qquad \Scal(\pi) := \int_\RRd \MCov(\pi_x,q_x) \,\mu(dx).
\end{equation}
where
\begin{equation}
\label{eq:MCov}
\MCov(\pi_x,q_x) := \sup_{p\in\Pi(\pi_x,q_x)} \int_{\RRd\times\RRd} m\cdot y \,p(dm,dy)
\end{equation}
denotes the maximum covariance between the measures $\pi_x$ and $q_x$; see Proposition \ref{prop:equiv-optimizations}.

In particular, in Section~\ref{sec:optimality}, we introduce the notion of a $q_x$-Bass martingale and prove that it is an optimizer of both \eqref{eq:mainPb} and \eqref{eq:mainPb2}. In contrast to the approach of Tschiderer and much of the related literature \cite{Ts24, BaScTs23}, our proof is direct and does not rely on duality theory.

In Section~\ref{sec:duality}, we establish a strong duality result for \eqref{eq:mainPb} and \eqref{eq:mainPb2}. More specifically, under suitable assumptions, we show that \eqref{eq:mainPb2} admits the following dual problem:
\[
    WT_{\Scal}^{q_\cdot}(\mu,\nu) = \inf_{\substack{\psi\in L^1(\nu)\\ \psi\ \mathrm{proper\ l.s.c.\ convex}}} \left(\int_\RRd \psi \, d\nu - \int_\RRd \psi^{C_{q_\cdot}} \, d\mu \right),
\]
where $\psi^{C_{q_\cdot}}$ denotes the $C_{q_\cdot}$-transform of the convex function $\psi$; see Definition~\ref{def:C_q-transform}.

Moreover, the equivalence between \eqref{eq:mainPb} and \eqref{eq:mainPb2} yields the following strong duality result for \eqref{eq:mainPb}:
\[
    WT_{\Ical}^{q_\cdot}(\mu,\nu) = \hspace{-4mm}\sup_{\substack{\psi\in L^1(\nu)\\ \psi\ \mathrm{proper}\\\mathrm{l.s.c.\ convex}}} 2 \left(\int_{\RRd \times \RRd}\hspace{-2mm} \psi^{C_{q_\cdot}}(x) + \frac{|y|^2}2  \,q_x(dy)  \, d\mu(x) - \int_\RRd \hspace{-2mm} \psi(x)-\frac{|x|^2}{2} \, d\nu(x)\right).
\]
We then introduce a generalization of the Martingale Sinkhorn algorithm developed in \cite{HaJoLoObPa26} and show that, whenever its iterates are well defined, it satisfies the same dual descent properties. Under additional compactness assumptions, we establish the convergence of the algorithm and, in particular, prove the existence of a $q_x$-Bass martingale between two given distributions; see Theorem~\ref{thm:q_x-sinkhorn-convergence}.
These results provide the foundation for the forthcoming work~\cite{Ma26_volMimicking}, which studies a dynamic martingale optimal transport problem aimed at calibrating martingale models to observed option prices while keeping their dynamics close to those of a prescribed reference model.

Finally, in Section~\ref{sec:1D}, we specialize the previous results to dimension one. Following \cite{AcMaPa23}, we derive a simplified representation of $q_x$-Bass martingales, characterize their existence through a fixed-point equation (Proposition \ref{prop:qx-bass-characterization}), and reformulate Algorithm~\ref{alg:q-sinkhorn} as a fixed-point iteration for a suitable operator $\Gcal_{q_\cdot}$ (Definitions \ref{def:G-qdot} and \ref{def:G-qx-iteration}).

In independent work, Beiglböck--Hasenbichler--Pammer \cite{BeHaPa26} introduce a closely related generalization of $q$-Bass martingales, named $R$-Bass martingales. While the $q_x$-Bass framework developed here allows for multidimensional marginals and reference kernels on $\RR^d$ indexed by $x\in\RR^d$, their construction considers one-dimensional marginal constraints and reference kernels on $\RR$ indexed by an augmented state space $\RR\times E$, where $E$ is a Polish space. In dimension one, the underlying weak martingale transport problems coincide when the reference kernels are independent of the $E$-coordinate, including the case where $E$ is a singleton. The two works establish complementary theoretical results under different assumptions; see \cite{BeHaPa26} for further details.

\paragraph{Notation.}
Let $\Prob(\RRd)$ denote the space of probability measures on $\RRd$. For $p\in[1,\infty)$, we write $\Prob_p(\RRd)$ for the collection of probability measures with finite $p$-th moment, whereas $\Prob_\infty(\RRd)$ denotes those with compact support. We use $\gamma^d$ for the standard Gaussian distribution on $\RRd$ and $\lambda^d$ for Lebesgue measure on $\RRd$. We denote by $\lambda$ the Lebesgue measure on $\RR$. Given $\xi \in \Prob(\RRd)$, the notation $\xi \ll \lambda^d$ means that $\xi$ is absolutely continuous with respect to $\lambda^d$, in which case its density is denoted by $\rho_\xi$. The support of $\xi$ is denoted by $\supp(\xi)$. If $\xi$ is a one-dimensional distribution, its cumulative distribution function is denoted by $F_\xi$, and its left-continuous quantile function is defined by $Q_\xi:(0,1)\to\RR$, $Q_\xi(u):=\inf\{y\in\RR:F_\xi(y)\geq u\}$. We extend it to $[0,1]$ by setting $Q_\xi(0):=\inf(\supp(\xi))$ and $Q_\xi(1):=\sup(\supp(\xi))$. Whenever $\xi\in\Prob_1(\RRd)$, its mean is denoted by $\mean(\xi):=\int_{\RRd} y\,\xi(dy)$.

For $\mu,\nu \in \Prob(\RRd)$, we denote by $\Pi(\mu,\nu)$ the set of probability measures on $\RRd\times\RRd$ whose first and second marginals are $\mu$ and $\nu$, respectively. Its elements are called couplings of $\mu$ and $\nu$. For $p \in [1,\infty)$ and $\xi,\zeta \in \Prob(\RRd)$, their $p$-Wasserstein distance is defined by
\[
    \Wcal_p(\xi,\zeta):=\inf_{\pi\in\Pi(\xi,\zeta)}\left(\int_{\RRd \times \RRd}|x-y|^p \, \pi(dx,dy)\right)^{1/p} \in [0,\infty].
\]
Given $\pi \in \Pi(\mu,\nu)$, we write $\pi(dx,dy)=\mu(dx)\pi_x(dy)$ for its disintegration with respect to the first marginal. We denote by $\Mart(\mu,\nu)$ the set of couplings $\pi \in \Pi(\mu,\nu)$ satisfying $\mean(\pi_x)=x$ for $\mu$-a.e.\ $x$, and refer to its elements as martingale couplings from $\mu$ to $\nu$.

Let $T:\RRd\to\RRd$ be measurable and $\xi\in\Prob(\RRd)$. The push-forward of $\xi$ under $T$ is denoted by $T_\#\xi$ and defined by $T_\#\xi(A):=\xi(T^{-1}(A))$ for every Borel set $A \subseteq \RRd$. For $\xi, \zeta \in \Prob(\RRd)$, their convolution is the probability measure $\xi \ast \zeta$ given by $(\xi\ast\zeta)(A):=\int_{\RRd\times\RRd}1_A(x+y)\,\xi(dx)\zeta(dy)$. For measurable functions $f$ and $g$, we set $f\ast g(x):=\int_{\RRd} f(x-y)g(y)\,dy$. Similarly, for a measurable function $f$ and $\xi\in\Prob(\RRd)$, we define $\xi\ast f(x):=\int_{\RRd} f(x-y)\,\xi(dy)$.

For a measure $\xi$ and $p\in[1,\infty]$, we denote by $L^p(\xi)$ the usual space of $p$-integrable functions with respect to $\xi$ (essentially bounded functions when $p=\infty$), with norm $\|\cdot\|_{L^p(\xi)}$. We write $L^0(A)$ for measurable real-valued functions on $A$, with representatives understood when evaluating or composing them. For $A\subseteq\RRd$, $L^\infty(A)$ denotes the space of $\lambda^d$-essentially bounded functions. The standard $L^\infty$ norm is denoted by $\|\cdot\|_\infty$. We use the shorthand $L^\infty(a,b):=L^\infty((a,b))$, with analogous conventions for $(a,b]$, $[a,b)$, and $[a,b]$. We denote by $\operatorname{co}(A)$ the convex hull of a set $A\subseteq\RRd$.

For a function $v:\RRd\to (-\infty, \infty]$, its convex conjugate is denoted by $v^*(y):=\sup_{x\in\RRd}{x\cdot y-v(x)}$.

\section{The Optimality of $q_x$-Bass Martingales}\label{sec:optimality}

This section aims to generalize the notion of $q$-Bass martingale introduced by Tschiderer in \cite{Ts24} and to show that its generalization is the optimizer of the weak optimal transport problems \eqref{eq:mainPb} and \eqref{eq:mainPb2}. These problems respectively generalize
\begin{equation}
\label{eq:Ts24_WOT}
    \inf_{\pi \in \Mart(\mu, \nu)} \int_\RRd \Wcal_2(\pi_x, q)^2\, \mu(dx) \qquad\text{and}\qquad \sup_{\pi \in \Mart(\mu, \nu)} \int_\RRd \MCov(\pi_x, q)\, \mu(dx),
\end{equation}
which were originally studied in \cite{Ts24}.
In particular, Tschiderer showed that, whenever it exists and under suitable integrability conditions on $\mu$, $\nu$, and $q$, a $q$-Bass martingale from $\mu$ to $\nu$ is a martingale coupling $\pi^* \in \Mart(\mu, \nu)$ that attains, respectively, the minimum and the maximum in the two problems in \eqref{eq:Ts24_WOT}.
The $q$-Bass construction can be outlined as follows. One starts with a suitable auxiliary random variable $A$ with law $\alpha\in\Prob(\RRd)$ and adds an independent noise $Q$ with distribution $q$. The sum $A+Q$ is then transported to the terminal distribution $\nu$ by the map $\nabla v$, where $v:\RRd\to\RR$ is convex. The initial state $M_0$ of the martingale is obtained by taking the conditional expectation of the terminal random variable with respect to $A$, namely
\[
    M_0=T_0(A):=\EE[\nabla v(A+Q)| A]=\int_\RRd \nabla v(A+z) \,q(dz).
\]

Thus, the initial marginal is given by $(T_0)_\#\alpha$, and the auxiliary distribution $\alpha$ must therefore be chosen so that $(T_0)_\#\alpha=\mu$. When this condition holds, $\alpha$ is called the Bass distribution, or equivalently the fixed-point distribution, for the $q$-Bass martingale from $\mu$ to $\nu$.

We now turn to the case considered in this paper, where the reference measure $q$ is replaced by a probability kernel $(q_x)_{x\in\supp(\mu)}$. In this setting, the reference noise can no longer be independent of the initial random variable. Instead, one considers
\[
    M_0\sim\mu, \qquad Q| M_0=x\sim q_x.
\]

As in the state-independent case, the terminal random variable is obtained by applying the gradient of a convex function. However, an additional Borel function $a:\RRd\to\RRd$ is now used to specify a random translation of the state-dependent noise $Q$:
\[ 
    M_1:=\nabla v(a(M_0)+Q).
\]

The function $a$ must be chosen, together with $v$, so that the resulting pair $(M_0,M_1)$ is a martingale. Equivalently, it must satisfy
\[
    x=\int_\RRd \nabla v\bigl(a(x)+z\bigr)\,q_x(dz) \qquad\text{for $\mu$-a.e.\ } x.
\]

The auxiliary Bass distribution is the law of $a(M_0)$, namely $\alpha=a_\#\mu$. Unlike in the state-independent setting, the relation between $\mu$ and $\alpha$ is therefore no longer required to be induced by the gradient of a convex function. Instead, it is encoded by the Borel function $a$. This motivates the following definition.

\begin{definition}[$q_x$-Bass martingale]
\label{def:qx-bass}
	Let $\mu,\nu\in\Prob_1(\RRd)$ and $(q_x)_{x\in\supp(\mu)} \subseteq \Prob(\RRd)$. We say that $\pi^*\in\Mart(\mu,\nu)$ is a $q_x$-Bass martingale from $\mu$ to $\nu$ if there exist a proper l.s.c. convex function $v:\RRd\to (-\infty, \infty]$ and a Borel function $a:\RRd\to\RRd$ such that, for $\mu$-a.e. $x$,
	\[
		v(a(x)+\cdot)\in L^1(q_x), \qquad \nabla v(a(x)+\cdot) \text{ exists } q_x\text{-a.s.},
	\]
	and, whenever $M_0 \sim \mu$ and $\Law(Q| M_0=x)=q_x$, one has
    \[
        \Law\bigl(M_0,\nabla v(a(M_0)+Q)\bigr)=\pi^*.
    \]
\end{definition}

More specifically, any $q_x$-Bass martingale from $\mu$ to $\nu$ makes the following diagram commute. Notice, in particular, that the random variable $Q$ depends on $M_0$, since its conditional distribution must satisfy $Q|M_0=x\sim q_x$. In the one-dimensional setting, a convenient way to construct $Q$ is to set $Q=f(M_0,U)$, where $f(x,u):=Q_{q_x}(u)$, $Q_{q_x}$ denotes the quantile function of $q_x$, and $U\sim\text{Unif}_{([0,1])}$ is independent of $M_0$.
\begin{figure}[H]
\centering
\begin{tikzcd}[row sep=0.7in, column sep=1.4in]
  M_0 \sim \mu
  \arrow[r, "\text{\normalsize $q_x$-Bass martingale}"]
  \arrow[d, "\text{\normalsize $a$}"']
  \arrow[dr, dashed, "\text{\normalsize $Q$ depends on $M_0$}"{sloped}]
  &
  M_1 \sim \nu
  \\
  a(M_0) \sim \alpha
  \arrow[r]
  &
  a(M_0)+Q
  \arrow[u, "\text{\normalsize $\nabla v$}"]
\end{tikzcd}
\caption{Diagram for the $q_x$-Bass martingale from $\mu$ to $\nu$}
\label{fig:diagram}
\end{figure}
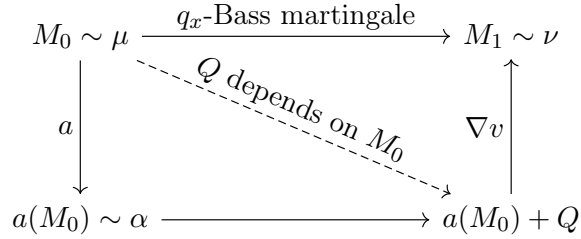

Any measure $\alpha$ making this diagram in Figure \ref{fig:diagram} commute will be referred to as a fixed-point distribution, or equivalently as a Bass distribution, in analogy with the standard $q$-Bass setting. Its existence necessarily requires $\mu$ to be dominated by $\nu$ in convex order; see \cite{St65}. Another structural condition that will play an important role in the subsequent sections is the irreducibility of the pair $(\mu,\nu)$.

\begin{definition}[Convex order and irreducibility]
\label{def:cvx-order-irreducibility}
Let $\mu,\nu\in\Prob_1(\RRd)$. We say that $\mu$ is dominated by $\nu$ in convex order, and write $\mu\preceq_c\nu$, if
\[
    \int_\RRd \varphi \, d\mu\leq\int_\RRd \varphi \, d\nu
\]
for every convex function $\varphi:\RRd \to \RR$ for which both integrals are well defined. The pair $(\mu,\nu)$ is called irreducible if $\mu\preceq_c\nu$ and, for every pair of Borel sets $A,B\subseteq\RRd$ with $\mu(A)>0$ and $\nu(B)>0$, there exists $\pi\in\Mart(\mu,\nu)$ such that
\begin{equation}
\label{eq:irreducible_coupling}
    \pi(A\times B)>0.
\end{equation}
\end{definition}

\begin{theorem}[Optimality of $q_x$-Bass martingales]
\label{thm:q_x-bass-optimality}
Let $p_1, p_2 \in [1, \infty]$ be Hölder conjugates, let $\mu,\nu\in\Prob_{p_1}(\RRd)$, and let $(q_x)_{x\in\supp(\mu)} \subseteq \Prob(\RRd)$. Suppose that there exists a $q_x$-Bass martingale $\pi^*\in\Mart(\mu,\nu)$ from $\mu$ to $\nu$ with Bass measure $\alpha \in \Prob_{p_2}(\RRd)$. Then the following statements hold:
\begin{enumerate}[label=(\roman*)]
    \item \label{it:optimality1} For any $\pi \in \Mart(\mu, \nu)$,
    \[
        \Ical(\pi^*)\leq\Ical(\pi).
    \]
    In particular, $WT_{\Ical}(\mu, \nu)$ is finite if and only if $\Ical(\pi^*)<\infty$, and, whenever this holds, $\pi^*$ is an optimizer of \eqref{eq:mainPb}.
    
    \item \label{it:optimality2} Assume that $p_1, p_2 < \infty$ and $\int_{\RRd\times\RRd}|z|^{p_2}\,q_x(dz)\mu(dx)<\infty$. Then
    \[
        \Scal(\pi^*) \geq \Scal(\pi)
    \]
    for any $\pi \in \Mart(\mu, \nu)$. In particular, $\pi^*$ is an optimizer of \eqref{eq:mainPb2}.
\end{enumerate}
\end{theorem}

\begin{remark}
    Under the assumptions of Theorem \ref{thm:q_x-bass-optimality}\ref{it:optimality2} $WT_{\Scal}(\mu, \nu)$ is finite. Indeed, if $p_1, p_2 \in [1, \infty)$ are Hölder conjugates and 
    \[
        M := \int_{\RRd\times\RRd} |z|^{p_2} q_x(dz)\mu(dx) < \infty.
    \]
    then Young's inequality yields
    \[
        \left| \int_\RRd \MCov(\pi_x, q_x) \, \mu(dx) \right| \leq 
        \frac{1}{p_1} \int_{\RRd\times\RRd} |m|^{p_1} \pi_x(dm)\mu(dx) + \frac{M}{p_2} < \infty,
    \]
    for all $\pi \in \Mart(\mu, \nu)$. Since the right-hand side is independent of $\pi$, the bound is uniform over $\pi \in \Mart(\mu, \nu)$.
\end{remark}

\begin{proof}
We begin by proving \ref{it:optimality1}. Let $\pi\in\Mart(\mu,\nu)$. If $\Ical(\pi)=+\infty$, there is nothing to prove. Hence, assume that $\Ical(\pi)<+\infty$. Let $(\pi_x)_{x \in \supp(\mu)}$ be a Borel disintegration of $\pi$ with respect to $\mu$. By \cite[Corollary 5.22]{Vi09}, we may choose a Borel kernel $(\eta_x)_{x \in \supp(\mu)} \subseteq \Prob(\RRd \times \RRd)$ such that
\[
    \eta_x \in \Pi(\pi_x,q_x) \qquad \text{and} \qquad \int_{\RRd\times\RRd}|m-z|^2 \,\eta_x(dm,dz) = \Wcal_2^2(\pi_x,q_x), \qquad \text{for $\mu$-a.e. $x$}.
\]

We now define a probability measure on $\RRd\times\RRd\times\RRd$ by $\PP(dx,dm,dz):=\mu(dx)\eta_x(dm,dz)$. Taking $M_0,M_1,Q$ to be the three coordinate projections, we have $M_0\sim\mu$ and $(M_1,Q)\mid M_0=x\sim\eta_x$. Then
\[
\Law(M_1\mid M_0=x)=\pi_x, \qquad \Law(Q\mid M_0=x)=q_x, \qquad \Ical(\pi)=\EE\left[|M_1-Q|^2\right].
\]
Moreover, since $\pi\in\Mart(\mu,\nu)$, $\Law(M_1)=\nu$ and $\EE[M_1\mid M_0]=M_0$.

Now let $\pi^*$ be a $q_x$-Bass martingale from $\mu$ to $\nu$ represented as in Definition \ref{def:qx-bass}, and let $\alpha \in \Prob_{p_2}(\RRd)$ denote its Bass measure. Thus, there exist a proper l.s.c. convex function $v:\RRd\to (-\infty, \infty]$ and a Borel map $a:\RRd\to\RRd$ such that
\[
    M_1^*:=\nabla v(a(M_0)+Q), \quad \EE[M_1^*\mid M_0]=M_0, \qquad (M_0,M_1^*)\sim\pi^*, \qquad a(M_0) \sim \alpha.
\]

For $\mu$-a.e. $x$, define $T_x(z):=\nabla u_x(z)$ and  $u_x(z):=v(a(x)+z)$. Since the map $z \mapsto u_x(z)$ is convex, $(\mathrm{id},T_x)_\#q_x$ is optimal between $q_x$ and $\pi_x^*:=(T_x)_\#q_x$ for the quadratic cost. By symmetry of the cost, $(T_x,\mathrm{id})_\#q_x$ is optimal between $\pi_x^*$ and $q_x$. Therefore,
\[
\Wcal_2^2(\pi_x^*,q_x) = \int_\RRd |T_x(z)-z|^2 \,q_x(dz)
\]
for $\mu$-a.e. $x$, with values in $[0,+\infty]$. Therefore, $\Ical(\pi^*) = \EE\left[|M_1^*-Q|^2\right]$.

Now, set $D = |M_1 - Q|^2 - |M_1^* - Q|^2$. Since $\Ical(\pi)=\EE[|M_1 - Q|^2] < \infty$ and $D^+ \leq |M_1 - Q|^2$, we have $D^+ \in L^1$. We next introduce
\[
\widetilde D = |M_1 - Q - a(M_0)|^2 - |M_1^* - Q - a(M_0)|^2.
\]
Expanding the squares gives $D = \widetilde D + 2a(M_0) \cdot (M_1-M_1^*)$. Since $a(M_0) \sim \alpha \in \Prob_{p_2}(\RRd)$ and $M_1, M_1^* \sim \nu \in \Prob_{p_1}(\RRd)$, Hölder's inequality yields
\[
    \EE[|2a(M_0) \cdot (M_1-M_1^*)|] < \infty.
\]
Hence, $\widetilde D^+ \in L^1$. Moreover, using the tower property and the martingale property of $(M_0, M_1)$ and $(M_0, M_1^*)$, $\EE[2a(M_0)\cdot(M_1-M_1^*)] = \EE[2a(M_0)\cdot \EE [M_1-M_1^*\mid M_0]] = 0$. Therefore,
\[
\Ical(\pi) - \Ical(\pi^*) = \EE[|M_1 - Q - a(M_0)|^2 - |M_1^* - Q- a(M_0)|^2].
\]

Recall now that Fenchel--Young gives
\[
v(Q+a(M_0))+v^*(M_1^*) = (Q+a(M_0)) \cdot M_1^*,
\]
whereas Fenchel--Young inequality yields $v(Q+a(M_0))+v^*(M_1) \geq (Q+a(M_0))\cdot M_1$.
Thus, we can conclude that
\[
\widetilde D \geq X-Y, \qquad X:=|M_1|^2-2v^*(M_1),\, Y:= |M_1^*|^2-2v^*(M_1^*).
\]

Since $\Law(M_1)=\Law(M_1^*)$, the CDFs of $X$ and $Y$ agree, and
\begin{align*}
    \EE[(X-Y)^+] & = \int_\RR \PP(Y \leq t < X) \, dt = \int_\RR F_Y(t) - \PP(X \leq t, Y \leq t) \, dt\\
    & = \int_\RR F_X(t) - \PP(X \leq t, Y \leq t) \, dt = \int_\RR \PP(X \leq t < Y)\, dt = \EE[(X-Y)^-].
\end{align*}
Since $0\leq(X-Y)^+ \leq \widetilde D^+ \in L^1$,
\[
\EE[(X-Y)^-] = \EE[(X-Y)^+]<\infty.
\]
Therefore, $X-Y$ is integrable, and $\EE[(X-Y)^-] = \EE[(X-Y)^+]$ implies that $\EE[X-Y]=0$. Thus,
\[
\Ical(\pi)-\Ical(\pi^*) \geq 0.
\]

The proof of \ref{it:optimality2} follows by the same argument. Fix $\pi \in \Mart(\mu, \nu)$. Since $WT_\Scal(\mu, \nu)$ is finite, $\Scal(\pi)<\infty$. Additionally, we may assume that $\Scal(\pi)>-\infty$. Otherwise, there is nothing to prove.

We repeat the same construction as in the proof of \ref{it:optimality1}, choosing this time a Borel kernel $(\eta_x)_{x \in \supp(\mu)} \subseteq \Prob(\RRd \times \RRd)$ such that
\[
    \eta_x \in \Pi(\pi_x,q_x) \qquad \text{and} \qquad \int_{\RRd\times\RRd}m \cdot z \,\eta_x(dm,dz) = \MCov(\pi_x,q_x), \qquad \text{for $\mu$-a.e. $x$}.
\]
In particular, \cite[Corollary 5.22]{Vi09} can also be applied in this case, since $c(m,z):=|m|^{p_1}/p_1+|z|^{p_2}/p_2-m\cdot z$ is continuous and non-negative, and minimizing this transport cost over $\Pi(\pi_x,q_x)$ is equivalent to maximizing $\int m\cdot z\,d\eta$.
Defining $\PP(dx,dm,dz):=\mu(dx)\eta_x(dm,dz)$ and taking $M_0,M_1,Q$ to be the three coordinate projections, we have
\[
\Scal(\pi) = \EE[M_1\cdot Q] \qquad \text{and} \qquad \Scal(\pi^*) = \EE[M_1^*\cdot Q],
\]
where $(M_0,M_1^*)$ is the $q_x$-Bass martingale constructed as in the proof of \ref{it:optimality1}.

We then set
\[
D = M_1\cdot Q-M_1^*\cdot Q, \qquad \widetilde D = M_1\cdot (Q+a(M_0))-M_1^*\cdot (Q+a(M_0)),
\]
so that $D = \widetilde D - a(M_0)\cdot(M_1-M_1^*)$. In particular, since $\Scal(\pi)>-\infty$ and $\Scal(\pi^*)<\infty$, we have
\[
D^- \leq (M_1\cdot Q)^-+(M_1^*\cdot Q)^+ \in L^1.
\]

As in the proof of \ref{it:optimality1}, $a(M_0)\cdot(M_1-M_1^*)$ is integrable and has expected value $0$. Therefore,
\[
\Scal(\pi)-\Scal(\pi^*) = \EE[D] = \EE[\widetilde D].
\]
Additionally, since $D^- \in L^1$, we have $\widetilde D^- \in L^1$. Finally, by the Fenchel--Young inequality,
\[
\widetilde D \leq v^*(M_1)-v^*(M_1^*).
\]
Since $\widetilde D^- \in L^1$ and $0 \leq (v^*(M_1)-v^*(M_1^*))^- \leq \widetilde D^-$, $(v^*(M_1)-v^*(M_1^*))^- \in L^1$. Therefore, using the same argument as above and the fact that $\Law(M_1)=\Law(M_1^*)$, we obtain
\[
(v^*(M_1)-v^*(M_1^*))^+ \in L^1 \qquad \text{and} \qquad \EE[v^*(M_1)-v^*(M_1^*)]=0.
\]
Thus, $\Scal(\pi)-\Scal(\pi^*) \leq 0$. 
\end{proof}

\begin{proposition}[Rescaling of the reference kernel]
Let $\mu,\nu\in\Prob_1(\RRd)$ and $(q_x)_{x\in\supp(\mu)}\subseteq\Prob(\RRd)$. For $s>0$, define $\widetilde q_x:=(s\,\mathrm{Id})_\# q_x$. If $\pi\in\Mart(\mu,\nu)$ is a $q_x$-Bass martingale, then $\pi$ is also a $\widetilde q_x$-Bass martingale.
\end{proposition}

\begin{proof}
Let $\pi\in \Mart(\mu,\nu)$ be a $q_x$-Bass martingale. Fix $s>0$ and let $v$, $a$, $M_0$, and $Q$ be as in Definition~\ref{def:qx-bass}, so that $\Law(M_0,\nabla v(a(M_0)+Q))=\pi$. Define $\widetilde v(z):=s v(z/s)$, $\widetilde a(x):=s a(x)$, and $\widetilde Q:=s Q$. Then $\Law(\widetilde Q\mid M_0=x)=(s\,\mathrm{Id})_\#q_x=\widetilde q_x$. Moreover, we have $\nabla\widetilde v(\widetilde a(M_0)+\widetilde Q)=\nabla v(a(M_0)+Q)$. The required integrability and differentiability conditions for $\widetilde v$ and $\widetilde a$ follow immediately from those of $v$ and $a$. Hence, by Definition~\ref{def:qx-bass}, $\Law(M_0,\nabla\widetilde v(\widetilde a(M_0)+\widetilde Q))=\pi$, and therefore $\pi$ is also a $\widetilde q_x$-Bass martingale.
\end{proof}

We conclude this section by recording a standard equivalence between \eqref{eq:mainPb} and \eqref{eq:mainPb2}, which will be used in the following section.

\begin{proposition} \label{prop:equiv-optimizations}
	Let $\mu,\nu\in\Prob_2(\RRd)$ and $(q_x)_{x\in\supp(\mu)} \subseteq \Prob(\RRd)$. Assume that
	\[
		\int_{\RRd\times\RRd} |y|^2 \,q_x(dy) \,\mu(dx) <\infty.
	\]
	Then $\pi^*\in\Mart(\mu,\nu)$ is an optimizer of $WT_{\Ical}^{q_\cdot}(\mu,\nu)$ if and only if it is an optimizer of $WT_{\Scal}^{q_\cdot}(\mu,\nu)$.
\end{proposition}

\begin{proof}
	For any $\pi\in\Mart(\mu,\nu)$, expanding the square gives
	\[
		\Ical(\pi) = \int_\RRd|x|^2 \,\nu(dx) + \int_\RRd \int_\RRd|y|^2 \,q_x(dy) \,\mu(dx) - 2\Scal(\pi).
	\]
	The first two terms are finite by assumption and do not depend on $\pi$.
	Thus minimizing $\Ical$ over $\Mart(\mu,\nu)$ is equivalent to maximizing $\Scal$ over $\Mart(\mu,\nu)$. Hence $\pi^*$ optimizes $WT_{\Ical}^{q_\cdot}(\mu,\nu)$ if and only if it optimizes $WT_{\Scal}^{q_\cdot}(\mu,\nu)$.
\end{proof}

\section{Strong Duality and Generalized Martingale Sinkhorn Algorithm} \label{sec:duality}
This section introduces the dual formulations of the weak optimal transport problems \eqref{eq:mainPb} and \eqref{eq:mainPb2}, which will serve as the basis for the numerical scheme studied in this section. The strong duality result for \eqref{eq:mainPb2} is obtained by adapting the arguments of \cite{BaBeScTs23, Ts24, HaJoLoObPa26} to the present state-dependent setting and is stated under the following assumptions.

\begin{assumption}[Assumptions for the duality result]
\label{ass:dual}
Let $p_1, p_2\in(1,\infty)$ be Hölder conjugates. Let $\mu,\nu\in\Prob_{p_1}(\RRd)$ and $(q_x)_{x\in\supp(\mu)}\subseteq\Prob_{p_2}(\RRd)$.
We assume that:
\begin{enumerate}[label=(A\arabic*)]
    \item \label{ass:dual1} the map $x\mapsto q_x$ is continuous with respect to the $p_2$-Wasserstein distance;
    \item \label{ass:dual2} if $p_1=p_2=2$, then $q_x\ll \lambda^d$ for $\mu$-a.e. $x\in\supp(\mu)$;
    \item \label{ass:dual3} 
    \[
        \int_{\RRd\times\RRd} |y|^{p_2} \, q_x(dy) \,\mu(dx) < \infty .
    \]
\end{enumerate}
\end{assumption}

To formulate the corresponding dual problems, we first introduce a state-dependent extension of the $C$-transform appearing in the dual formulation of the martingale Benamou--Brenier problem. Recall that, if $\gamma^d$ denotes the standard Gaussian distribution, the $C$-transform of a proper l.s.c.\ convex function $\psi$ is defined by
\[
\psi^C(x):=(\psi^* \ast \gamma^d)^*(x)=\sup_{a\in\RRd} \left(a\cdot x -\int_\RRd \psi^*(a+z) \,\gamma^d(dz)\right),
\]
see \cite{HaJoLoObPa26}, for instance. The $C_{q_\cdot}$-transform generalizes this construction by replacing the fixed Gaussian measure $\gamma^d$ in the integral with the state-dependent reference kernel $(q_x)_{x \in \supp(\mu)}$. In particular, if $q_x=\gamma^d$ for every $x\in\supp(\mu)$, then $\psi^{C_{q_\cdot}}(x)=\psi^C(x)$.

\begin{definition}[$C_{q_\cdot}$-transform]\label{def:C_q-transform}
Let $\Omega \subseteq \RRd$, $(q_x)_{x\in \Omega} \subseteq \Prob(\RRd)$, and let $\psi:\RRd\to(-\infty, \infty]$ be a proper l.s.c. convex function. We define the $C_{q_\cdot}$-transform of $\psi$ by
\[
    \psi^{C_{q_\cdot}}(x) := \sup_{a\in\RRd} \left ( a\cdot x - \int_\RRd \psi^*(a+z) \,q_x(dz)\right ),
\]
for every $x\in\Omega$ for which the right-hand side is well-defined.
\end{definition}

\begin{theorem}[Dual formulation of \eqref{eq:mainPb2}]
\label{thm:dual}
Let $p_1, p_2\in(1,\infty)$ be Hölder conjugates, $\mu,\nu\in\Prob_{p_1}(\RRd)$ and $(q_x)_{x\in\supp(\mu)} \subseteq \Prob_{p_2}(\RRd)$ and assume that $\mu\preceq_c\nu$. Under Assumption \ref{ass:dual}, the primal problem \eqref{eq:mainPb2} admits an optimizer and has finite value. Moreover,
\begin{equation}\label{eq:dual}
    WT_{\Scal}^{q_\cdot}(\mu,\nu) = \inf_{\substack{\psi\in L^1(\nu)\\ \psi\ \textrm{proper l.s.c.}\\\textrm{convex}}} \Ecal^{q_\cdot}(\psi), \qquad \Ecal^{q_\cdot}(\psi) := \int_\RRd \psi \, d\nu - \int_\RRd \psi^{C_{q_\cdot}} \, d\mu.
\end{equation}
Additionally, if $p_1=p_2=2$, the optimizer of $WT_{\Scal}^{q_\cdot}(\mu,\nu)$ is unique.
\end{theorem}

Since the proof is rather technical and follows by adapting arguments from \cite{Ts24, HaJoLoObPa26}, we postpone it to Appendix \ref{sec:appendix}. The dual formulation of \eqref{eq:mainPb} then follows directly from Theorem \ref{thm:dual} and the equivalence established in Proposition \ref{prop:equiv-optimizations}.

\begin{corollary}[Dual formulation of \eqref{eq:mainPb}]
\label{cor:dual_cor}
Let $\mu,\nu\in\Prob_2(\RRd)$ and $(q_x)_{x\in\supp(\mu)} \subseteq \Prob_2(\RRd)$ such that $\mu\preceq_c\nu$. Under Assumption \ref{ass:dual} with $p_1=p_2=2$, the primal problem \eqref{eq:mainPb} admits a unique optimizer and has finite value. Moreover,
\begin{equation}\label{eq:dual_cor}
    WT_{\Ical}^{q_\cdot}(\mu,\nu) = 2 \cdot \hspace{-4mm}  \sup_{\substack{\psi\in L^1(\nu)\\ \psi\ \mathrm{proper}\\\mathrm{l.s.c.\ convex}}}  \left(\int_{\RRd \times \RRd}\hspace{-2mm} \psi^{C_{q_\cdot}}(x) + \frac{|y|^2}2  \,q_x(dy)  \, d\mu(x) - \int_\RRd \hspace{-2mm} \psi(x)-\frac{|x|^2}{2} \, d\nu(x)\right).
\end{equation}
\end{corollary}

\begin{proof}
The proof follows immediately from Theorem \ref{thm:dual} and Proposition \ref{prop:equiv-optimizations}. Indeed, Assumption \ref{ass:dual} ensures that all the hypotheses of Proposition \ref{prop:equiv-optimizations} are satisfied.
\end{proof}

We next use \eqref{eq:dual} to introduce a generalized Martingale Sinkhorn algorithm for computing a $q_x$-Bass distribution. This extends the Martingale Sinkhorn algorithm of \cite{HaJoLoObPa26} from the classical Bass setting to the present state-dependent framework.

\begin{breakablealgorithm}
\caption{Generalized Martingale Sinkhorn algorithm}
\label{alg:q-sinkhorn}
\noindent Let $\psi_0 \in L^1(\nu)$ be a proper l.s.c. convex function.\\

\noindent For $k\geq 1$:
\begin{enumerate}[label=(\roman*)]
    \item \label{it:q-sinkhorn-1} Define a Borel map $a_k:\supp(\mu) \to \RRd$ such that
    \begin{equation}
        \label{eq:alg_eq1}
         a_k(x) \in \operatorname*{arg \,max}_{a\in\RRd} \left(a \cdot x - \int_\RRd \psi_{k-1}^*(a+z) \,q_x(dz)\right).\vspace{-2mm}
    \end{equation}
    so that $\psi_{k-1}^{C_{q_\cdot}}(x) = a_k(x)\cdot x - \int_\RRd \psi_{k-1}^*(a_k(x)+z) \,q_x(dz)$.

    \item \label{it:q-sinkhorn-2} Define
    \begin{equation}
        \label{eq:alg_eq2}
        \rho_k := (a_k(x)+z)_\# \left(q_x(dz)\mu(dx)\right).
    \end{equation}

    \item \label{it:q-sinkhorn-3} Choose a proper l.s.c. convex function $\psi_k$ such that $\psi_k^*$ is a Brenier potential from $\rho_k$ to $\nu$, i.e.
    \begin{equation}
        \label{eq:alg_eq3}
        \nu=(\nabla\psi_k^*)_\#\rho_k.
    \end{equation}
\end{enumerate}
\end{breakablealgorithm}

As in the original Martingale Sinkhorn algorithm, we start from a convex potential $\psi_0$ and alternate between two optimization problems. Given $\psi_{k-1}$, the optimization problem \eqref{eq:alg_eq1} determines the state-dependent translation $a_k(x)$ of the reference measures, for each $x\in\supp(\mu)$. Then, in step \ref{it:q-sinkhorn-2}, the map $a_k$ then determines the translated reference measure $\rho_k$. The optimization problem in step \ref{it:q-sinkhorn-3} updates the potential by solving the quadratic optimal transport problem from $\rho_k$ to the prescribed terminal distribution $\nu$.

Under very general assumptions, these operations need not always be well defined. For instance, the convolution of $\psi_{k-1}^*$ with $q_x$ may fail to be finite. The following algorithm is therefore understood whenever all the required operations are well defined.

Whenever the iterates of the Generalized Martingale Sinkhorn Algorithm are well defined, the same descent mechanism as in the original Martingale Sinkhorn algorithm shows that $(\psi_k)_{k\in\NN}$ is a descent sequence for the dual objective in \eqref{eq:dual}. Moreover, if no further decrease occurs at some iteration, then the dual potential, and consequently the corresponding Bass distribution, has been reached. These properties are summarized in the following theorem.

\begin{theorem}[Descent property]
\label{thm:q_x-sinkhorn-descent}
Let $\mu,\nu \in \Prob_1(\RRd)$ and $(q_x)_{x\in\supp(\mu)} \subseteq \Prob(\RRd)$ satisfying $q_x \ll \lambda^d$ for $\mu$-a.e. $x \in \supp(\mu)$. Let $(\psi_k)_{k \in\NN}$ denote the sequence of potentials generated by Algorithm \ref{alg:q-sinkhorn}, and let $\Ecal^{q_\cdot}$ denote the dual objective in \eqref{eq:dual}. Assume that all steps of Algorithm \ref{alg:q-sinkhorn} are well-defined and that, for every $k\in\NN$,
\[
\int_\RRd (|\psi_{k-1}|+|\psi_k|)\,d\nu + \int_\RRd (|\psi_{k-1}^{C_{q_\cdot}}|+|\psi_k^{C_{q_\cdot}}|)\,d\mu + \int_\RRd (|\psi_{k-1}^*|+|\psi_k^*|)\,d\rho_k + \int_\RRd |a_k(x)\cdot x|\,\mu(dx)<\infty.
\]
Then
\begin{equation}
\label{eq:descent_property}
    \Ecal^{q_\cdot}(\psi_k) \leq \Ecal^{q_\cdot}(\psi_{k-1}), \qquad k\in\NN.
\end{equation}

Moreover, fix $k\in\NN$ such that equality holds in \eqref{eq:descent_property}. If $\rho_k\in\Prob_{p_2}(\RRd)$ and $\nu\in\Prob_{p_1}(\RRd)$ for some Hölder conjugates $p_1,p_2\in[1,\infty]$ and, for $\mu$-a.e. $x$,
\[
\nabla_a\int_{\RRd}\psi_{k-1}^*(a+z)\,q_x(dz) = \int_{\RRd}\nabla\psi_{k-1}^*(a+z)\,q_x(dz)
\]
at $a=a_k(x)$, then $\psi_{k-1}^*$ generates a $q_x$-Bass martingale from $\mu$ to $\nu$
\[
    \pi^*(dx,dy) = \mu(dx) \, \left(\nabla\psi_{k-1}^* (a_k(x) + \cdot)\right)_\#q_x(dy).
\]
\end{theorem}

The proof of the theorem follows by a straightforward adaptation of the descent argument in \cite[Lemma 3.5]{HaJoLoObPa26} to accommodate the state-dependent reference kernel. For completeness, we provide the details in Appendix \ref{sec:appendix}.

\section{Convergence of Algorithm \ref{alg:q-sinkhorn} and Existence of $q_x$-Bass Martingales}
In this section, we show that, under the assumptions below, the sequence of potentials $(\psi_k)_{k \in \NN}$ generated by Algorithm \ref{alg:q-sinkhorn} is well-defined. Moreover, under the additional assumption that the pair $(\mu,\nu)$ is irreducible (see Definition \ref{def:cvx-order-irreducibility}), we show that $(\psi_k)_{k \in \NN}$ can be normalized so as to be subsequentially compact. This allows us to prove that every accumulation point is a potential generating a $q_x$-Bass martingale from $\mu$ to $\nu$. Our argument follows the strategy of \cite{HaJoLoObPa26}, adapted to the setting of state-dependent kernels.

\begin{assumption}[Assumptions for the compact setting]
\label{ass:q_x-sinkhorn-compact}
Let $\mu,\nu\in\Prob_\infty(\RRd)$ and $(q_x)_{x\in\supp(\mu)}\subseteq\Prob_1(\RRd)$. We assume that:
\begin{enumerate}[label=(C\arabic*)]
    \item \label{ass:compact1} the map $x\mapsto q_x$ is continuous with respect to the $1$-Wasserstein distance;
    \item \label{ass:compact2} $q_x\ll\lambda^d$ for every $x\in\supp(\mu)$;
    \item \label{ass:compact3} $\mean(q_x)=x$ for every $x\in\supp(\mu)$;
    \item \label{ass:compact4} $\supp(\mu)\subseteq\interior(\co(\supp(\nu)))$.
\end{enumerate}
\end{assumption}

\begin{proposition}[Well-posedness and compactness of the iterates]
\label{prop:q_x-sinkhorn-compactness}
Under Assumption~\ref{ass:q_x-sinkhorn-compact}, Algorithm~\ref{alg:q-sinkhorn} is well defined for every initial potential $\psi_0\in L^1(\nu)$ such that $\psi_0^*$ is Lipschitz.

If, moreover, $(\mu,\nu)$ is irreducible, let $(\psi_k)_{k\in\NN}$ and $(a_k)_{k\in\NN}$ denote, respectively, the sequences of potentials and Borel functions generated by Algorithm~\ref{alg:q-sinkhorn}, and set $\bar x:=\mean(\nu)$. For each $k\in\NN$, choose $p_k\in\partial\psi_k(\bar x)$ so that
\[
    \ell_k(y):=p_k\cdot(y-\bar x)+\psi_k(\bar x)
\]
is an affine minorant of $\psi_k$ and define $\widehat\psi_k:=\psi_k-\ell_k$ and $\widehat a_{k+1}:=a_{k+1}-p_k$. Then we have
\[
    \sup_{k\in\NN}\|\widehat a_{k+1}\|_{L^\infty(\mu)}<\infty.
\]
In particular, every subsequence admits a further subsequence $(k_n)_{n\in\NN}$ and a proper l.s.c.\ convex function $\widehat\psi$, finite on $\interior(\co(\supp(\nu)))$, such that
\begin{equation}
\label{eq:epiconvergence_conclusion}
    \widehat\psi_{k_n}\to\widehat\psi \quad \text{locally uniformly on } \interior(\co(\supp(\nu))), \qquad
    \widehat\psi_{k_n}^*\to\widehat\psi^* \quad \text{locally uniformly on } \RRd.
\end{equation}
\end{proposition}

We postpone the proof of Proposition~\ref{prop:q_x-sinkhorn-compactness} to Appendix~\ref{sec:appendix} and focus here on the main result of this section.

\begin{theorem}[Convergence and existence of a $q_x$-Bass martingale]
\label{thm:q_x-sinkhorn-convergence}
Under Assumption~\ref{ass:q_x-sinkhorn-compact}, assume that $(\mu,\nu)$ is irreducible.

Let $(\psi_k)_{k\in\NN}$ and $(a_k)_{k\in\NN}$ be, respectively, the sequences of potentials and Borel functions generated by Algorithm~\ref{alg:q-sinkhorn} from an initial potential $\psi_0\in L^1(\nu)$ such that $\psi_0^*$ is Lipschitz. Let $(\widehat\psi_k)_{k\in\NN}$ and $(\widehat a_{k+1})_{k\in\NN}$ denote the normalized sequences defined as in Proposition~\ref{prop:q_x-sinkhorn-compactness}.

Then every accumulation point $\widehat\psi^*$ of $(\widehat\psi_k)_{k\in\NN}$ generates a $q_x$-Bass martingale from $\mu$ to $\nu$.
\end{theorem}

\begin{proof}
By Proposition~\ref{prop:q_x-sinkhorn-compactness}, Algorithm~\ref{alg:q-sinkhorn} is well defined. Furthermore, combining \eqref{eq:estimate_dual_obj} with $\mu\preceq_c\nu$, we obtain
\[
    \Ecal^{q_\cdot}(\psi_k) \geq \int_\RRd \psi_k \, d\nu - \int_\RRd \psi_k \, d\mu + \int_\RRd |x|^2 \,\mu(dx) \geq \int_\RRd |x|^2 \,\mu(dx).
\]
Hence, the sequence $(\Ecal^{q_\cdot}(\psi_k))_{k \in \NN}$ is bounded from below. On the other hand, by Theorem~\ref{thm:q_x-sinkhorn-descent}, $(\Ecal^{q_\cdot}(\psi_k))_{k \in \NN}$ is decreasing. It therefore converges, and in particular,
\begin{equation}
\label{eq:q_x-energy-gap-vanishes}
    \Ecal^{q_\cdot}(\psi_k) - \Ecal^{q_\cdot}(\psi_{k+1}) \longrightarrow 0, \qquad \text{as } k \to \infty.
\end{equation}

By Proposition~\ref{prop:q_x-sinkhorn-compactness}, there exist a subsequence $(\widehat\psi_{k_n})_{n \in \NN}$ and a proper l.s.c. convex function $\widehat\psi$ such that \eqref{eq:epiconvergence_conclusion} holds. For $n \in \NN$, define
\begin{equation}
\label{eq:G_def_existence_convergence}
    G_n(x,a):=\int_\RRd \widehat\psi_{k_n}^*(a+z) \,q_x(dz), \qquad G(x,a):=\int_\RRd \widehat\psi^*(a+z) \,q_x(dz).
\end{equation}
Since $\supp(\mu)$ is compact and $x\mapsto q_x$ is continuous in $\Wcal_1$, the family $(q_x)_{x \in \supp(\mu)}$ is compact in $\Prob_1(\RRd)$ and therefore has uniformly integrable first moments. Hence, using the local uniform convergence of $\widehat\psi_{k_n}^*$ to $\widehat\psi^*$ and their common Lipschitz constant, we obtain
\begin{equation}
\label{eq:G_n_uniform_convergence}
    \sup_{\substack{x \in \supp(\mu)\\ a \in C}} |G_n(x,a)-G(x,a)| \longrightarrow 0
\end{equation}
for every compact set $C \subseteq \RRd$. Moreover, $|G_n(x,a)-G_n(y,b)| \leq L|a-b|+L\Wcal_1(q_x,q_y)$, and the same estimate holds for $G$, so that $G_n$ and $G$ are continuous on $\supp(\mu) \times \RRd$.

Proposition~\ref{prop:q_x-sinkhorn-compactness} also yields a compact set $K \subseteq \RRd$ such that
\[
    \widehat a_{k+1}(x) \in K, \qquad \text{for all } k \in \NN,\ x \in \supp(\mu).
\]
Since this does not in general yield a pointwise convergent subsequence of $(\widehat a_{k_n+1})_{n \in \NN}$, we consider instead
\[
    \eta_n:=(id,\widehat a_{k_n+1})_{\#}\mu.
\]
As $\supp(\eta_n)\subseteq\supp(\mu)\times K$, the sequence $(\eta_n)_{n \in \NN}$ is tight. Thus, up to a further subsequence, $(\eta_n)_{n \in \NN}$ weakly converges to $\eta$, for some $\eta \in \Prob_\infty(\RRd\times\RRd)$ with $\supp(\eta)\subseteq\supp(\mu)\times K$. Since the first marginal of $\eta_n$ is $\mu$ for every $n \in \NN$, the same holds for $\eta$, and we may write
\[
    \eta(dx,da)=\alpha_x(da)\mu(dx),
\]
where $\alpha_x\in\Prob(K)$ for $\mu$-a.e. $x$.

We claim that
\begin{equation}
\label{eq:alpha_x_optimality}
    \alpha_x\left(\operatorname*{argmax}_{a \in \RRd}\big(a\cdot x-G(x,a)\big)\right)=1, \qquad \mu\text{-a.e. }x \in \supp(\mu).
\end{equation}
Indeed, for any $n\in\NN$, set $m_n, m: \supp(\mu)\times \RRd \to \RR$
\[
    m_n(x,a):=a\cdot x-G_n(x,a), \qquad m(x,a):=a\cdot x-G(x,a),
\]
and $M_n, M: \supp(\mu) \to \RR$
\[
    M_n(x):=\sup_{a \in \RRd}m_n(x,a), \qquad M(x):=\sup_{a \in \RRd}m(x,a).
\]
By the optimality of $\widehat a_{k_n+1}$ and the fact that $\widehat a_{k_n+1}(x)\in K$,
\[
    M_n(x)=m_n(x,\widehat a_{k_n+1}(x))=\max_{a \in K}m_n(x,a).
\]
Hence, by \eqref{eq:G_n_uniform_convergence},
\[
    \sup_{x \in \supp(\mu)}|M_n(x)-M(x)| \longrightarrow 0, \qquad \text{as } n \to \infty.
\]
where, passing to the limit in $m_n(x,a)\leq M_n(x)$ for arbitrary $a \in \RRd$, we also obtain $M(x)=\max_{a \in K}m(x,a)$. Therefore,
\[
    \int_{\RRd\times\RRd}\big(M_n(x)-m_n(x,a)\big)\,\eta_n(dx,da)
    \longrightarrow
    \int_{\RRd\times\RRd}\big(M(x)-m(x,a)\big)\,\eta(dx,da),
\]
since $M_n-m_n\to M-m$ uniformly on $\supp(\mu)\times K$ and $M-m$ is continuous there. On the other hand,
\[
    \int_{\RRd\times\RRd}\big(M_n(x)-m_n(x,a)\big)\,\eta_n(dx,da)
    =\int_\RRd \big(M_n(x)-m_n(x,\widehat a_{k_n+1}(x))\big)\,\mu(dx)=0.
\]
Thus,
\[
    \int_\RRd\int_K \big(M(x)-m(x,a)\big)\,\alpha_x(da)\mu(dx)=0.
\]
Since $M(x)-m(x,a)\geq0$, we conclude that
\[
    M(x)=m(x,a), \qquad \alpha_x\text{-a.e. }a \in K,\quad \text{for }\mu\text{-a.e. }x \in \supp(\mu),
\]
which proves \eqref{eq:alpha_x_optimality}.

Now, for every $k\in\NN$, set
\[
    \rho_{k+1}:=(a_{k+1}(x)+z)_\#\bigl(q_x(dz)\mu(dx)\bigr), \quad \text{and}\quad \widehat\rho_{k+1}:=(\widehat a_{k+1}(x)+z)_\#\bigl(q_x(dz)\mu(dx)\bigr).
\]
Then, the weak convergence of $(\eta_n)_{n \in \NN}$ to $\eta$, together with the continuity of $x \mapsto q_x$ in $\Wcal_1$ and the uniform boundedness of $(\widehat a_{k_n+1})_{n \in \NN}$, yields
\[
    \widehat\rho_{k_n+1} \longrightarrow \rho:=(a+z)_\#\bigl(q_x(dz)\alpha_x(da)\mu(dx)\bigr), \qquad \text{in } \Wcal_1 \text{ as } n \to \infty.
\]
Moreover, since $q_x \ll \lambda^d$ for every $x \in \supp(\mu)$, we also have $\rho \ll \lambda^d$.

Furthermore, applying \eqref{eq:inequality3_dual_proof} with $k=k_n+1$ and using the optimality of $\psi_{k_n+1}$ in the dual formulation of $\MCov(\rho_{k_n+1},\nu)$ yields
\[
0 \leq \int_\RRd \psi_{k_n}^* \, d\rho_{k_n+1} + \int_\RRd \psi_{k_n} \, d\nu-\MCov(\rho_{k_n+1},\nu) \leq \Ecal^{q_\cdot}(\psi_{k_n})-\Ecal^{q_\cdot}(\psi_{k_n+1}).
\]
Then the invariance under the simultaneous affine normalization of $\psi_{k_n}$ and translation of $\rho_{k_n+1}$ yields
\begin{equation}
\label{eq:q_x-brenier-gap}
0 \leq \int_\RRd \widehat \psi_{k_n}^* \, d\widehat\rho_{k_n+1} + \int_\RRd \widehat\psi_{k_n} \, d\nu - \MCov(\widehat\rho_{k_n+1},\nu) \longrightarrow 0
\end{equation}
by \eqref{eq:q_x-energy-gap-vanishes}. Since $\widehat\rho_{k_n+1} \to \rho$ in $\Wcal_1$, $\nu$ is compactly supported and $\widehat\psi_{k_n}^* \to \widehat\psi^*$ locally uniformly with a common Lipschitz constant,
\[
    \int_\RRd \widehat\psi_{k_n}^* \, d\widehat\rho_{k_n+1} \longrightarrow \int_\RRd \widehat\psi^* \, d\rho,\qquad \text{as } n \to \infty,
\]
by the dominated convergence theorem.
Additionally, $\MCov(\widehat\rho_{k_n+1},\nu)$ converges to $\MCov(\rho,\nu)$ as $n \to \infty$, since $\nu \in \Prob_\infty(\RRd)$ yields $|\MCov(\widehat\rho_{k_n+1},\nu) - \MCov(\rho,\nu)|\leq \sup_{y \in \supp(\nu)} |y| \Wcal_1(\widehat\rho_{k_n+1}, \rho)$.

Moreover, $\widehat\psi_{k_n}\geq0$ and $\widehat\psi_{k_n}$ epi-converges to $\widehat\psi$, so Fatou's lemma gives
\[
    \int_\RRd \widehat\psi \, d\nu \leq \liminf_{n\to\infty} \int_\RRd \widehat\psi_{k_n} \, d\nu.
\]
Therefore, passing to the limit in \eqref{eq:q_x-brenier-gap} yields
\[
    \int_\RRd \widehat\psi^* \, d\rho+\int_\RRd \widehat\psi \, d\nu\leq\MCov(\rho,\nu).
\]
In particular, this shows that $\widehat \psi$ attains the minimum in \eqref{eq:dual_MCov}. Hence, equality holds, and since $\rho\ll\lambda^d$,
\begin{equation}
\label{eq:final_transport_cond}
    \nu=(\nabla\widehat\psi^*)_\#\rho.
\end{equation}

It remains to verify that $\widehat\psi^*$ generates a $q_x$-Bass martingale from $\mu$ to $\nu$. Let $(M_0, A, Q)$ be such that
\[
\Law(M_0, A, Q)= \mu(dx)\alpha_x(da)q_x(dz)
\]
and set $M_1:= \nabla \widehat \psi^*(A+Q)$.  By \eqref{eq:final_transport_cond}, $M_1 \sim \nu$. 

Moreover, since $q_x \ll \lambda^d$ and $\widehat \psi^*$ is Lipschitz for $\mu$-a.e. $x$, $G(x , \cdot)$ is differentiable and satisfies
\[
\nabla_a G(x, a)= \int_\RRd \nabla \widehat \psi^*(a+z) q_x(dz).
\]
On the other hand, by \eqref{eq:alpha_x_optimality}, $\alpha_x$ is concentrated on the maximizers of $m(x,\cdot)$ for $\mu$-a.e. $x$. Hence, the first-order optimality condition yields $x = \nabla_a G(x,a)$ for $\alpha_x$-a.e. $a \in K$, for $\mu$-a.e. $x$. Therefore, we have 
\begin{equation}
\label{eq:martingale_property_existence}
\EE[\EE[M_1|M_0, A]|M_0]=\EE[M_0|M_0]=M_0.
\end{equation}
In particular, $\pi^* = \Law(M_0, M_1) \in \Mart(\mu, \nu)$.

Therefore, it remains to show that $\widehat \psi^*$ generates a $q_x$-Bass martingale from $\mu$ to $\nu$. To this end, we first prove that $\widehat\psi$ minimizes $\Ecal^{q_\cdot}$. Definition \ref{def:C_q-transform} and \eqref{eq:alpha_x_optimality} give
\[
\widehat\psi^{C_{q_\cdot}}(x)+\int_\RRd \widehat\psi^*(a+z) \, q_x(dz)=a\cdot x,\qquad \alpha_x\text{-a.e. }a.
\]
Since $\Law(M_0,A,Q)=\mu(dx)\alpha_x(da)q_x(dz)$ and $M_1\sim\nu$, integrating the previous identity with respect to $\mu$ yields
\[
\int_\RRd \widehat\psi^{C_{q_\cdot}} \, d\mu = \EE\left[A\cdot M_0-\widehat\psi^*(A+Q)\right],
\]
and therefore
\[
\Ecal^{q_\cdot}(\widehat\psi) = \EE \left[ \widehat \psi(M_1) + \widehat\psi^*(A+Q) -A \cdot M_0\right].
\]
Since $M_1 = \nabla \widehat \psi^*(A+Q)$, Fenchel--Young's equality gives $\widehat \psi(M_1) + \widehat \psi^*(A+Q) = (A+Q) \cdot M_1$. Hence, using \eqref{eq:martingale_property_existence} and the tower property,
\begin{equation}
\label{eq:connection_dual_objective}
\Ecal^{q_\cdot}(\widehat\psi) = \EE\left[(A+Q)\cdot M_1-A\cdot M_0\right] = \EE[Q\cdot M_1].
\end{equation}

Now, denote by $(\pi^*_x)_{x \in \supp(\mu)}$ a Borel disintegration of $\pi^*$. Conditionally on $M_0=x$, $\Law(M_1, Q) \in \Pi(\pi^*_x, q_x)$.
Therefore, by \eqref{eq:connection_dual_objective}, we have
\[
\Ecal^{q_\cdot}(\widehat\psi) \leq \int_\RRd \MCov(\pi^*_x, q_x) \mu(dx).
\]

Conversely, let $\varphi \in L^1(\nu)$ be proper l.s.c. convex and fix $a \in K$ and $\theta \in \Pi(\pi_x^*, q_x)$. By Fenchel--Young's inequality,
\[
y\cdot(a+z)\leq \varphi(y)+\varphi^*(a+z).
\]
Integrating with respect to $\theta$ and using that $\mean(\pi_x^*)=x$, we obtain
\[
\int_{\RRd\times\RRd} y\cdot z\,\theta(dy,dz) \leq \int_\RRd \varphi(y) \, \pi_x^*(dy)-a\cdot x+\int_\RRd \varphi^*(a+z)\,q_x(dz).
\]
Taking the supremum over $\theta \in \Pi(\pi_x^*, q_x)$, the infimum over $a \in \RRd$, and then integrating with respect to $\mu$ yields
\[
\int_\RRd \MCov(q_x,\pi_x^*)\,\mu(dx) \leq \int_\RRd \varphi\,d\nu-\int_\RRd \varphi^{C_{q_\cdot}}\,d\mu = \Ecal^{q_\cdot}(\varphi).
\]
Thus $\Ecal^{q_\cdot}(\widehat \psi) \leq \Ecal^{q_\cdot}(\varphi)$, for every admissible $\varphi$. So $\widehat \psi$ is a minimizer of $\Ecal^{q_\cdot}$.

Since $\widehat \psi \in L^1(\nu)$ and $\widehat \psi^*$ is Lipschitz, Proposition \ref{prop:q_x-sinkhorn-compactness} allows us to apply Algorithm \ref{alg:q-sinkhorn} starting from the potential $\widehat \psi$. If $\widetilde \psi$ denotes the potential obtained after one iteration, then Theorem \ref{thm:q_x-sinkhorn-descent} and the minimality of $\widehat \psi$ for the dual objective $\Ecal^{q_\cdot}$ imply $\Ecal^{q_\cdot}(\widetilde \psi) = \Ecal^{q_\cdot}(\widehat \psi)$. Thus, by the equality case in Theorem \ref{thm:q_x-sinkhorn-descent}, $\widehat \psi^*$ generates a $q_x$-Bass martingale from $\mu$ to $\nu$.
\end{proof}

\section{The One-Dimensional Setting}\label{sec:1D}
In this section, following the approach of \cite{AcMaPa23}, we introduce a simplified representation for $q_x$-Bass martingales in dimension one. This setting is particularly relevant for applications. In mathematical finance, for instance, one may seek a martingale coupling between prescribed marginal distributions, inferred from European option prices through the Breeden--Litzenberger formula \cite{BrLi78}, while imposing specific requirements on its dynamics. Since this formula determines the one-dimensional distribution of the underlying asset at each maturity, such calibration problems naturally fit into a one-dimensional martingale transport framework.

As in the case of $q$-Bass martingales, we characterize the existence of a $q_x$-Bass martingale between two prescribed marginals through an equation. Besides providing a convenient theoretical description, this characterization yields a simplified version of Algorithm~\ref{alg:q-sinkhorn} in dimension one. Following \cite{AcMaPa23}, we now define, for each $A\in L^0(0,1)$, the map $S_A^{\mu,\nu,q_\cdot}:(0,1)\times\RR\to\RR$, which will be used to formulate both the equation and the algorithm.

\begin{definition}
\label{def:S-qdot}
Let $\mu,\nu\in\Prob_1(\RR)$ and $(q_x)_{x\in\supp(\mu)} \subseteq \Prob(\RR)$ such that $q_x\ll\lambda$ for $\mu\text{-a.e. }x$. For $A \in L^0(0,1)$, we define
\[
    S_A^{\mu, \nu, q_\cdot}(u,x):=\int_\RR Q_\nu\left(\int_0^1 F_{q_{Q_\mu(v)}}(x+z-A(v)) \, dv\right)\rho_{q_{Q_\mu(u)}}(z) \, dz.
\]
\end{definition}

To motivate the introduction of the map $S_A^{\mu,\nu,q_\cdot}$, consider the construction in Definition \ref{def:qx-bass} and parametrize the initial state as $M_0=Q_\mu(U)$, where $U\sim \text{Unif}_{[0,1]}$. The measurable function $A:(0,1)\to\RR$ represents the state-dependent translation, in the sense that $A(u)=a(Q_\mu(u))$. If, conditionally on $U=u$, the noise $Q$ has distribution $q_{Q_\mu(u)}$, then the function
\[
y\longmapsto \int_0^1 F_{q_{Q_\mu(v)}}\bigl(y-A(v)\bigr) \, dv
\]
is the CDF of $A(U)+Q$. Consequently, its monotone rearrangement onto $\nu$ is given by
\begin{equation}
\label{eq:transport_map_1D}
T_A(y):=Q_\nu\left(\int_0^1 F_{q_{Q_\mu(v)}}\bigl(y-A(v)\bigr) \, dv\right),
\end{equation}
which plays in dimension one the role of the gradient $\nabla v$ in Definition \ref{def:qx-bass}. The corresponding terminal random variable is therefore $M_1=T_A(A(U)+Q)$, and its conditional mean at the initial state $Q_\mu(u)$ is precisely
\[
\EE[M_1| U=u]=S_A^{\mu,\nu,q_\cdot}(u,A(u)).
\]
Thus, the martingale condition is equivalent to the identity \eqref{eq:fixed-point-eq}.

\begin{proposition}[One-dimensional characterization]
\label{prop:qx-bass-characterization}
Let $p_1,p_2\in[1,\infty]$ be Hölder conjugates, let $\mu,\nu\in\Prob_{p_1}(\RR)$, and let $(q_x)_{x\in\supp(\mu)} \subseteq \Prob_{p_2}(\RR)$ such that $q_x\ll\lambda$ for $\mu\text{-a.e. }x$. Then there exists a $q_x$-Bass martingale from $\mu$ to $\nu$ if and only if there exists $a \in L^0(\RR)$ such that the map $A \in L^0(0,1)$ defined by $A = a \circ Q_\mu$ satisfies
\begin{equation}
\label{eq:fixed-point-eq}
    S_A^{\mu,\nu,q_\cdot}(u,A(u))=Q_\mu(u),\qquad \text{for a.e. }u\in(0,1).
\end{equation}
\end{proposition}

\begin{remark}[The factorization condition]
\label{rmk:factorization}
The factorization $A= a \circ Q_\mu$ in Proposition~\ref{prop:qx-bass-characterization} is not required for optimality. Indeed, any solution $A$ of \eqref{eq:fixed-point-eq} makes the $q_x$-Bass diagram in Figure \ref{fig:diagram} commute and yields optimality for \eqref{eq:mainPb} and \eqref{eq:mainPb2}, whenever the assumptions of Theorem~\ref{thm:q_x-bass-optimality} are satisfied. Rather, the role of the factorization condition is to identify a canonical parametrization in which the state-dependent translation is a measurable function of the initial state. Without this condition, as we showed in the proof of Theorem \ref{thm:q_x-sinkhorn-convergence}, the law of the Bass distribution associated with a $q_x$-Bass martingale does not even need to be induced by a measurable map as it can also be defined via a suitable kernel $(\alpha_x)_{x \in \supp(\mu)}$.
\end{remark}

\begin{remark}[Failure of quantile and CDF representations]
\label{rmk:failure_quantile_CDF}
In contrast with the state-independent setting, a solution $A\in L^0(0,1)$ of \eqref{eq:fixed-point-eq} need not be monotone and therefore need not be a quantile function. In the case of $q$-Bass martingales, one may choose $A$ as the quantile function of the Bass distribution, and the factorization condition $A = a \circ Q_\mu$ amounts to requiring that the map $q\star T$ appearing in \cite[Definition 1.1]{AcMa26_semidiscrete} is increasing $\alpha$-a.e., where $\alpha$ denotes the associated Bass distribution.

In the state-dependent setting, $A(u)$ records not only the distribution of the translations but also how each translation is coupled with the reference measure $q_{Q_\mu(u)}$ associated with the initial state $Q_\mu(u)$. Rearranging the values of $A$ may therefore change the distribution of the sum $A(U)+Q$, even if the law of $A$ remains unchanged. Thus, $A$ cannot in general be replaced by its monotone rearrangement, and the fixed-point formulation of Conze and Henry-Labordère \cite{CoHe21}, which is expressed in terms of CDFs, does not extend directly to $q_x$-Bass martingales.
\end{remark}

\begin{proof}
Let $a\in L^0(\RR)$ be such that $A=a\circ Q_\mu$ satisfies \eqref{eq:fixed-point-eq}, and let $U\sim\text{Unif}_{[0,1]}$. We set $M_0:=Q_\mu(U)$, $Q| U=u\sim q_{Q_\mu(u)}$, and define
\[
M_1:=Q_\nu\left(\int_0^1 F_{q_{Q_\mu(v)}}(A(U)+Q-A(v)) \, dv\right).
\]
By construction, $M_1\sim\nu$. Moreover, $\EE[M_1| U=u]=S_A^{\mu,\nu,q_\cdot}(u,A(u))=Q_\mu(u)$ for a.e. $u\in(0,1)$. Thus $\EE[M_1| U]=Q_\mu(U)=M_0$. Since $M_0$ is a function of $U$, $\EE[M_1| M_0]=\EE[\EE[M_1| U]| M_0]=M_0$.

Finally, let $v$ be a convex potential for the monotone transport $T_A$. By convexity, for some $C>0$,
\[
|v(y)|\le C(1+|y|)(1+|T_A(y)|)
\]
for $\Law(A(U)+Q)$-a.e. $y$. Since $\nu\in\Prob_{p_1}(\RR)$ and $q_x\in\Prob_{p_2}(\RR)$, H\"older's inequality yields $v(a(x)+\cdot)\in L^1(q_x)$ for $\mu$-a.e. $x$. Hence, the law of $(M_0,M_1)$ is a $q_x$-Bass martingale from $\mu$ to $\nu$.

Conversely, suppose that a $q_x$-Bass martingale exists. Then there are a proper l.s.c. convex function $v:\RR\to (-\infty, \infty]$ and a Borel map $a:\RR\to\RR$ such that
\[
    M_0\sim\mu,\qquad Q| M_0=x\sim q_x,\qquad M_1=\nabla v(a(M_0)+Q)\sim \nu, \quad \text{and} \quad M_0=\EE[M_1| M_0].
\]
Let $U\sim\text{Unif}_{[0,1]}$ and write $M_0=Q_\mu(U)$. Define $A(u):=a(Q_\mu(u))$. Since in dimension one the monotone transport is unique, we may identify
\[
    \nabla v(x)=Q_\nu\left(\int_0^1 F_{q_{Q_\mu(w)}}(x-A(w)) \, dw\right)
\]
for $\Law(A(U)+Q)$-a.e. $x$. Therefore, for a.e. $u\in(0,1)$,
\[
    Q_\mu(u)=\EE[M_1| U=u]=S_A^{\mu,\nu, q_\cdot}(u,A(u)).
\]
\end{proof}

We next reformulate Algorithm~\ref{alg:q-sinkhorn} as a fixed-point iteration.

\begin{definition}[The fixed-point operator $\Gcal_{q_\cdot}$]
\label{def:G-qdot}
Let $\mu,\nu\in\Prob_\infty(\RR)$ and $(q_x)_{x\in\supp(\mu)} \subseteq \Prob_1(\RR)$ such that $q_x \ll \lambda$ $\mu$-a.e. and $\supp(\mu) \subseteq \interior(\co(\supp(\nu)))$.
For any measurable function $A:(0,1)\to\RR$, we define
\[
\bigl(\Gcal_{q_\cdot}A\bigr)(u):=\bigl(S_A^{q_\cdot}(u,\cdot)\bigr)^{-1}\bigl(Q_\mu(u)\bigr),\qquad u\in(0,1),
\]
where
\[
\bigl(S_A^{q_\cdot}(u,\cdot)\bigr)^{-1}(y):=\inf\bigl\{x\in\RR:S_A^{q_\cdot}(u,x)\geq y\bigr\},\qquad y\in\RR,
\]
denotes the left-continuous generalized inverse of the non-decreasing map $x\mapsto S_A^{q_\cdot}(u,x)$.
\end{definition}

\begin{remark}
\label{rem:G-qx-left-inverse}
Under the assumptions of Definition \ref{def:G-qdot}, if
\[
    A = \Gcal_{q_\cdot} A,
\]
then
\[
    S_A^{q_\cdot}(u,A(u))=Q_\mu(u), \qquad \text{for a.e. }u\in(0,1),
\]
with $A=a\circ Q_\mu$ for some $a\in L^0(\RR)$. Therefore, by Proposition \ref{prop:qx-bass-characterization}, the existence of a fixed point of $\Gcal_{q_\cdot}$ implies the existence of a $q_x$-Bass martingale from $\mu$ to $\nu$.
\end{remark}

\begin{definition}[Fixed-point iteration]
\label{def:G-qx-iteration}
Under the assumptions of Definition \ref{def:G-qdot}, let $A_0:(0,1)\to\RR$ be measurable. The fixed-point iteration associated with $\Gcal_{q_\cdot}$ is defined by
\[
A_{k+1}:=\Gcal_{q_\cdot}A_k,\qquad k\ge0.
\]
\end{definition}

\begin{remark}[Equivalence between the fixed-point iteration and Algorithm \ref{alg:q-sinkhorn}]
\label{rmk:algorithmic-1D-equivalence}
The fixed-point iteration in Definition \ref{def:G-qx-iteration} is equivalent to Algorithm \ref{alg:q-sinkhorn}. Indeed, suppose that at iteration $k \geq 1$ we are given the measurable function $A_k$. By Definition \ref{def:G-qx-iteration}, this is equivalent to having a measurable map $a_k:\RR\to\RR$ such that $A_k=a_k\circ Q_\mu$. Thus, the fixed-point iteration can be viewed as starting from step \ref{it:q-sinkhorn-2} of Algorithm \ref{alg:q-sinkhorn}. At this step, one determines the law $\rho_k$, which is equivalent to computing the integral
\[
\int_0^1 F_{q_{Q_\mu(v)}}(A_k(u)+z-A_k(v))\,dv
\]
appearing in the definition of $S_{A_k}^{q_\cdot}$. Then, applying $Q_\nu$ to this function yields the transport map from $\rho_k$ to $\nu$, corresponding to step \ref{it:q-sinkhorn-3} of Algorithm \ref{alg:q-sinkhorn}.

Finally, $A_{k+1}$ is obtained by applying $\bigl(S_{A_k}^{q_\cdot}(u,\cdot)\bigr)^{-1}$ to $Q_\mu(u)$. This is exactly the relation used in step \ref{it:q-sinkhorn-1} of the next iteration to define the measurable map $a_{k+1}$. More precisely, for $x=Q_\mu(u)$, step \ref{it:q-sinkhorn-1} amounts to determining $A_{k+1}(u):=a_{k+1}(Q_\mu(u))\in\RR$ solving
\[
\max_{a\in\RR}\left(aQ_\mu(u)-(\psi_k^*\ast q_{Q_\mu(u)})(a)\right).
\]
The corresponding first-order condition is
\[
Q_\mu(u)=\frac{d}{da}\left(\psi_k^*\ast q_{Q_\mu(u)}\right)(A_{k+1}(u))=S_{A_k}^{q_\cdot}(u,A_{k+1}(u)),
\]
and therefore
\[
A_{k+1}(u)=\bigl(S_{A_k}^{q_\cdot}(u,\cdot)\bigr)^{-1}(Q_\mu(u))=\bigl(\Gcal_{q_\cdot}A_k\bigr)(u),
\]
which is precisely the fixed-point iteration of Definition \ref{def:G-qx-iteration}.
\end{remark}

In view of its equivalence with Algorithm \ref{alg:q-sinkhorn}, the fixed-point scheme inherits the conclusions of Theorem \ref{thm:q_x-sinkhorn-convergence}, under the assumptions of that theorem.

\begin{remark}[Linear convergence]
Under stronger regularity assumptions, the fixed-point iteration associated with $\Gcal_{q_\cdot}$ converges linearly in the $L^\infty$ norm. The proof follows the arguments developed in \cite{AcMaPa23} for the Bass-LV model and adapted in \cite{AcMa26_general} to the $q$-Bass setting. The same argument extends to the present state-dependent framework, provided that $\nu$ has a density bounded away from zero on $\co(\supp(\nu))$, each $q_x$ admits a strictly positive density $\rho_x\in L^2(\RR)$, and the maps $x\mapsto\rho_x$ and $x\mapsto q_x$ satisfy the corresponding continuity assumptions. In particular, $\Gcal_{q_\cdot}$ is non-expansive on $L^\infty(0,1)$ and the fixed-point iteration converges linearly to the solution set.
\end{remark}

\appendix
\section{Postponed Proofs and Auxiliary Results}
\label{sec:appendix}

\begin{proof}[Proof of Theorem \ref{thm:dual}]
Let $C_{p_1}(\RRd)$ be the set of continuous functions with $p_1$-growth $\psi : \RRd \rightarrow \RR$ such that
\[
	C_{p_1}(\RRd) = \left\{\psi \in C(\RRd; \RR) \colon \ell+\frac{|\cdot|^{p_1}}{p_1} \leq \psi(\cdot)\leq a + k |\cdot|^{p_1}, \quad\text{for some } a, k, \ell \in \RR \right\}.
\]
For $\psi\in C_{p_1}(\RRd)$, define
\[
    \varphi^\psi(x):=\inf_{\substack{p\in\Prob_{p_1}(\RRd)\\\mean(p)=x}}\left(\int_\RRd \psi \, dp-\MCov(p,q_x)\right),\qquad x\in\supp(\mu), 
\]
 and set $m_1(x):=\mean(q_x)$ and $m_{p_2}(x):=\int_\RRd |z|^{p_2} \,q_x(dz)$. By Assumption \ref{ass:dual}, $m_{p_2}\in L^1(\mu)$, and by Hölder's and Jensen's inequalities, respectively,
 \begin{align*}
     \int_\RRd |x\cdot m_1(x)| \,\mu(dx) &\leq \left(\int_\RRd |x|^{p_1} \,\mu(dx)\right)^{1/p_1}\left(\int_\RRd |m_1(x)|^{p_2} \,\mu(dx)\right)^{1/p_2}\\
     &\leq \left(\int_\RRd |x|^{p_1} \,\mu(dx)\right)^{1/p_1}\left(\int_\RRd m_{p_2}(x) \,\mu(dx)\right)^{1/p_2} <\infty.
 \end{align*}
For $x \in \supp(\mu)$ and $p \in \Prob_{p_1}(\RRd)$, we define
\[
C(x,p) := \begin{cases} -\MCov(p,q_x) + \frac{1}{p_1} \int_\RRd |y|^{p_1} \, p(dy), & \text{if } \mean(p)=x,\\ +\infty, & \text{otherwise.} \end{cases}
\]
and set
\[
\widetilde C(x,p) := C(x,p) + \frac{1}{p_2} m_{p_2}(x).
\]
Then, whenever $\mean(p)=x$,
\[
\widetilde C(x,p)=\inf_{\pi \in \Pi(p,q_x)} \int_{\RRd\times\RRd} \left( \frac{|y|^{p_1}}{p_1} + \frac{|z|^{p_2}}{p_2} - y \cdot z\right) \pi(dy,dz)\geq 0,
\]
where the last inequality follows from Young's inequality. Since the cost $c(y,z):=\frac{|y|^{p_1}}{p_1}+\frac{|z|^{p_2}}{p_2}-y\cdot z$ is continuous and non-negative, the stability of optimal transport plans, see \cite[Theorem~5.20]{Vi09}, together with the convergence of the corresponding $p_1$- and $p_2$-moments, implies that the associated optimal transport value is continuous with respect to the $\Wcal_{p_1}\times\Wcal_{p_2}$ topology. Hence, since $x\mapsto q_x$ is $\Wcal_{p_2}$-continuous and the barycenter map is continuous for the $\Wcal_{p_1}$ topology, $\widetilde C$ is jointly l.s.c. Moreover, $\widetilde C\geq0$ and $\widetilde C(x,\cdot)$ is convex. Thus, \cite[Theorem~1.2]{BaBePa18} applies and yields the existence of an optimizer for
\[
\inf_{\pi \in \Pi(\mu,\nu)} \int_\RRd \widetilde C(x,\pi_x) \, \mu(dx),
\]
and therefore for $WT_\Scal^{q_\cdot}(\mu,\nu)$.

Let $C_{\mathrm b,p_1}(\RRd)$ denote the set of continuous functions with  $p_1$-growth from above and bounded from below by a constant. Then, by \cite[Theorem~1.3]{BaBePa18},
\[
\inf_{\pi \in \Pi(\mu,\nu)} \int_\RRd \widetilde C(x,\pi_x) \, \mu(dx) = \sup_{\widetilde\psi \in C_{\mathrm b,p_1}(\RRd)} \left\{\int_\RRd \inf_{p \in \Prob_{p_1}(\RRd)} \left(\int_\RRd \widetilde\psi \, dp + \widetilde C(x,p)\right) \mu(dx) - \int_\RRd \widetilde\psi \, d\nu\right\}.
\]
Setting $\psi(y) := \widetilde\psi(y) + \frac{1}{p_1} |y|^{p_1}$, we have $\psi \in C_{p_1}(\RRd)$ and
\[
\inf_{p \in \Prob_{p_1}(\RRd)} \left(\int_\RRd \widetilde\psi \, dp + \widetilde C(x,p)\right) = \varphi^\psi(x) + \frac1{p_2} m_{p_2}(x).
\]
Since
\[
\inf_{\pi \in \Pi(\mu,\nu)} \int_\RRd \widetilde C(x,\pi_x) \, \mu(dx) = \frac1{p_1} \int_\RRd |y|^{p_1} \, \nu(dy) + \frac1{p_2} \int_\RRd m_{p_2}(x) \, \mu(dx) - WT_{\Scal}^{q_\cdot}(\mu,\nu),
\]
we can conclude that
\begin{equation}
\label{eq:first_dual}
WT_{\Scal}^{q_\cdot}(\mu,\nu) = \inf_{\psi \in C_{p_1}(\RRd)} \left(\int_\RRd \psi \, d\nu - \int_\RRd \varphi^\psi \, d\mu\right).
\end{equation}

Additionally, if $p_1=p_2=2$, the assumption $q_x\ll\lambda^d$, for $\mu$-a.e. $x$, and Brenier's theorem imply that $\widetilde C(x,p) = \frac12\Wcal_2^2(p,q_x)$ on the affine constraint $\mean(p)=x$, and $p\mapsto\widetilde C(x,p)$ is strictly convex. Therefore, if $p_1=p_2=2$, the primal optimizer is unique. Indeed, if two optimizers had different conditional laws on a set of positive $\mu$-measure, their midpoint would still be admissible and would strictly improve the corresponding minimization of $\int C(x,\pi_x) \,\mu(dx)$.

We next reduce the first dual problem \eqref{eq:first_dual} to convex test functions. Let $\psi\in C_{p_1}(\RRd)$ and let conv$(\psi)$ be its convex hull. The measurable selection used in \cite[Proposition 2.3]{Ts24} depends only on the representation of conv$(\psi)$ and not on the reference measure. The same argument applies with $C_2(\RRd)$ and $\Prob_2(\RRd)$ replaced by $C_{p_1}(\RRd)$ and $\Prob_{p_1}(\RRd)$, respectively. In particular, the $p_1$-growth condition ensures that the measures obtained through the measurable selection have finite $p_1$-moment. The only estimate involving the reference measure is the monotonicity of $\MCov(\cdot,q)$ with respect to convex order, which remains valid for measures in $\Prob_{p_1}(\RRd)$ and $\Prob_{p_2}(\RRd)$, and it applies pointwise with $q=q_x$. More specifically, if $p_x\preceq_c\bar p_x$, then
\[
    \MCov(p_x,q_x)\leq\MCov(\bar p_x,q_x).
\]
Hence
\[
    WT_{\Scal}^{q_\cdot}(\mu,\nu)=\inf_{\substack{\psi\in C_{p_1}(\RRd)\\ \psi\ \mathrm{convex}}}\left(\int_\RRd \psi \, d\nu-\int_\RRd \varphi^\psi \, d\mu\right).
\]

We now identify $\varphi^\psi$ with the $C_{q_\cdot}$-transform. Fix a convex $\psi\in C_{p_1}(\RRd)$ and $x\in\supp(\mu)$, and define, for $\bar x\in\RRd$,
\[
    \Phi_x^\psi(\bar x):=\inf_{\substack{p\in\Prob_{p_1}(\RRd)\\ \mean(p)=\bar x}}\left(\int_\RRd \psi\,dp-\MCov(p,q_x)\right).
\]
Thus, $\varphi^\psi(x)=\Phi_x^\psi(x)$. Our goal is to identify $\Phi_x^\psi$ explicitly.

For a fixed reference measure $q$, \cite[Lemma~2.6]{Ts24} computes the corresponding optimization problem without the barycenter constraint, namely
\[
    \sup_{p\in\Prob_2(\RRd)}\left(\MCov(p,q)-\int_\RRd \eta\,dp\right)=\int_\RRd \eta^*\,dq
\]
for every proper convex function $\eta$. Its proof extends to the present setting by replacing $\Prob_2(\RRd)$ with $\Prob_{p_1}(\RRd)$ and $L^2(q)$ with $L^{p_1}(q)$, with Hölder's inequality replacing the corresponding quadratic integrability estimates. Applying this adapted identity with $q=q_x$ and $\psi_a(y):=\psi(y)-a\cdot y$, we obtain
\[
    (\Phi_x^\psi)^*(a)=\sup_{\bar x\in\RRd}\left(a\cdot\bar x-\Phi_x^\psi(\bar x)\right) =\int_\RRd \psi^*(a+z)\,q_x(dz).
\]

It remains to recover $\Phi_x^\psi$ from its convex conjugate. In particular, the argument in \cite[Proposition~2.7]{Ts24} establishing the applicability of the Fenchel--Moreau theorem can be adapted to the present setting using Young's inequality, which gives
\[
    \Phi_x^\psi(\bar x)\geq \ell-\frac{1}{p_2}m_{p_2}(x)>-\infty,
\]
while choosing $p=\delta_{\bar x}$ shows that $\Phi_x^\psi(\bar x)<\infty$. Since $\Phi_x^\psi$ is convex and finite on $\RRd$, it is continuous, and hence Fenchel--Moreau yields
\[
    \Phi_x^\psi(\bar x)=\sup_{a\in\RRd}\left(a\cdot\bar x-\int_\RRd \psi^*(a+z)\,q_x(dz)\right).
\]
Evaluating at $\bar x=x$, we conclude that
\[
    \varphi^\psi(x)=\sup_{a\in\RRd}\left(a\cdot x-\int_\RRd \psi^*(a+z)\,q_x(dz)\right)=\psi^{C_{q_\cdot}}(x).
\]

Now, recall that $\psi^*$ is continuous and has at most $p_2$-growth. Hence, by the $\Wcal_{p_2}$-continuity of $x\mapsto q_x$, for every $a\in\RRd$ the map
\[
    x\mapsto\int_\RRd \psi^*(a+z) \,q_x(dz)
\]
is continuous on $\supp(\mu)$, while for every $x\in\supp(\mu)$ it is continuous in $a$. Therefore, $x\mapsto\psi^{C_{q_\cdot}}(x)$ is Borel as the supremum over $a\in\QQ^d$ of Borel functions. Thus
\[
    WT_{\Scal}^{q_\cdot}(\mu,\nu)=\inf_{\substack{\psi\in C_{p_1}(\RRd)\\ \psi\ \mathrm{convex}}}\left(\int_\RRd \psi \, d\nu-\int_\RRd \psi^{C_{q_\cdot}} \, d\mu\right).
\]

It remains to pass from convex functions in $C_{p_1}(\RRd)$ to proper l.s.c. convex functions in $L^1(\nu)$. Fix an arbitrary $\pi(dx,dy)=\mu(dx)\pi_x(dy)\in\Mart(\mu,\nu)$ and define the relaxed dual functional
\[
    \Ecal^{q_\cdot}(\psi):=\int_\RRd \left[\int_\RRd \psi(y) \,\pi_x(dy)-\psi^{C_{q_\cdot}}(x)\right]\mu(dx),
\]
for proper l.s.c. convex $\psi$ satisfying $\mu(\operatorname{dom}\psi)=1$. The value of $\Ecal^{q_\cdot}(\psi)$ is independent of the chosen martingale transport, because the first term only depends on the second marginal $\nu$. Moreover, since $\psi$ is proper l.s.c. and convex, it admits an affine minorant, which controls its negative part.

We first check that $\Ecal^{q_\cdot}(\psi)$ is well-defined with an integrable lower bound. Fix $a\in\RRd$, $x\in\operatorname{dom}\psi$ and consider the proper convex function $\psi_a(y):=\psi(y)-a\cdot y$. Its convex conjugate is $\psi_a^*(z)=\psi^*(a+z)$. Hence, by Fenchel--Young's inequality,
\begin{equation} \label{eq:fenchel-young-inequality}
	    y\cdot z\leq\psi_a(y)+\psi_a^*(z)=\psi(y)-a\cdot y+\psi^*(a+z).
\end{equation}
Taking $y=x$ in \eqref{eq:fenchel-young-inequality} and integrating with respect to $q_x(dz)$ gives
\[
    a\cdot x-\int_\RRd \psi^*(a+z) \,q_x(dz)\leq\psi(x)-x\cdot m_1(x).
\]
Since $\mu(\operatorname{dom}\psi)=1$, taking the supremum over $a\in\RRd$ yields
\[
    \psi^{C_{q_\cdot}}(x)\leq\psi(x)-x\cdot m_1(x),\qquad \mu\text{-a.e. }x .
\]
By Jensen's inequality,
\[
    \int_\RRd \psi(y) \,\pi_x(dy)-\psi^{C_{q_\cdot}}(x)\geq x\cdot m_1(x),\qquad \mu\text{-a.e. }x,
\]
and the right-hand side is $\mu$-integrable by the estimate at the beginning of the proof. Thus $\Ecal^{q_\cdot}(\psi)$ is well-defined as an element of $\RR \cup \{+\infty\}$.

We claim that
\[
    WT_{\Scal}^{q_\cdot}(\mu,\nu)=\inf_{\substack{\psi\ \mathrm{proper\ l.s.c.\ convex}\\ \mu(\operatorname{dom}\psi)=1}}\Ecal^{q_\cdot}(\psi).
\]
Since any element of $C_{p_1}(\RRd)$ is proper l.s.c. and finite,
\[
    WT_{\Scal}^{q_\cdot}(\mu,\nu)= \inf_{\substack{\psi\in C_{p_1}(\RRd)\\ \psi\ \mathrm{convex}}}\left(\int_\RRd \psi \, d\nu-\int_\RRd \psi^{C_{q_\cdot}} \, d\mu\right) \geq \inf_{\substack{\psi\ \mathrm{proper\ l.s.c.\ convex}\\ \mu(\operatorname{dom}\psi)=1}}\Ecal^{q_\cdot}(\psi).
\]
Conversely, let $\psi$ be proper l.s.c. convex with $\mu(\operatorname{dom}\psi)=1$, let $x\in\operatorname{dom}\psi$, and let $\rho\in\Prob_{p_1}(\RRd)$ satisfy $\mean(\rho)=x$. Integrating \eqref{eq:fenchel-young-inequality} with respect to an arbitrary coupling $\eta\in\Pi(\rho,q_x)$, and using $\mean(\rho)=x$, we obtain
\[
    \int_{\RRd\times\RRd}y\cdot z \,\eta(dy,dz)\leq\int_\RRd \psi(y) \,\rho(dy)-a\cdot x+\int_\RRd \psi^*(a+z) \,q_x(dz),
\]
where the left-hand side is finite by Hölder's inequality. Taking the supremum over all couplings $\eta\in\Pi(\rho,q_x)$ gives
\[
    \MCov(\rho,q_x)\leq\int_\RRd \psi(y) \,\rho(dy)-a\cdot x+\int_\RRd \psi^*(a+z) \,q_x(dz).
\]
Taking the infimum over $a\in\RRd$ yields
\[
    \MCov(\rho,q_x)\leq\int_\RRd \psi(y) \,\rho(dy)-\psi^{C_{q_\cdot}}(x).
\]
Applying this inequality with $\rho=\pi_x$ and integrating with respect to $\mu(dx)$ gives
\[
    \int_\RRd \MCov(\pi_x,q_x) \,\mu(dx)\leq\Ecal^{q_\cdot}(\psi).
\]
Since the right-hand side is independent of the chosen martingale transport, taking the supremum over $\pi\in\Mart(\mu,\nu)$ gives $WT_{\Scal}^{q_\cdot}(\mu,\nu)\leq\Ecal^{q_\cdot}(\psi)$. Taking the infimum over $\psi$ proves the claim.

Now let $\psi\in L^1(\nu)$ be proper l.s.c. convex. Since $\mu\preceq_c\nu$, the convex-order relation implies that $\psi$ is $\mu$-integrable and that $\mu(\operatorname{dom}\psi)=1$. More explicitly, choosing an affine minorant $\ell$ of $\psi$, the nonnegative convex function $\psi-\ell$ satisfies
\[
    \int_\RRd(\psi-\ell) \, d\mu\leq\int_\RRd(\psi-\ell) \, d\nu<\infty,
\]
and therefore $\psi\in L^1(\mu)$. Moreover, the estimate above shows that $\int_\RRd \psi^{C_{q_\cdot}} \, d\mu$ is well-defined as an extended integral. For such $\psi$, the relaxed functional satisfies
\[
    \Ecal^{q_\cdot}(\psi)=\int_\RRd \psi \, d\nu-\int_\RRd \psi^{C_{q_\cdot}} \, d\mu.
\]
Consequently, using the inclusions
\begin{align*}
	\left\{\psi\in C_{p_1}(\RRd):\psi\ \mathrm{convex}\right\} & \subseteq\left\{\psi\in L^1(\nu):\psi\ \mathrm{proper\ l.s.c.\ convex}\right\}\\
	& \subseteq\left\{\psi\ \mathrm{proper\ l.s.c.\ convex}:\mu(\operatorname{dom}\psi)=1\right\},
\end{align*}
we obtain
\begin{align*}
    WT_{\Scal}^{q_\cdot}(\mu,\nu)&=\inf_{\substack{\psi\ \mathrm{proper\ l.s.c.\ convex}\\ \mu(\operatorname{dom}\psi)=1}}\Ecal^{q_\cdot}(\psi)\leq\inf_{\substack{\psi\in L^1(\nu)\\ \psi\ \mathrm{proper\ l.s.c.\ convex}}}\left(\int_\RRd \psi \, d\nu-\int_\RRd \psi^{C_{q_\cdot}} \, d\mu\right)\\
    &\leq\inf_{\substack{\psi\in C_{p_1}(\RRd)\\ \psi\ \mathrm{convex}}}\left(\int_\RRd \psi \, d\nu-\int_\RRd \psi^{C_{q_\cdot}} \, d\mu\right)=WT_{\Scal}^{q_\cdot}(\mu,\nu).
\end{align*}
Hence
\[
    WT_{\Scal}^{q_\cdot}(\mu,\nu)=\inf_{\substack{\psi\in L^1(\nu)\\ \psi\ \mathrm{proper\ l.s.c.}\\ \mathrm{convex}}}\left(\int_\RRd \psi \, d\nu-\int_\RRd \psi^{C_{q_\cdot}} \, d\mu\right).
\]
\end{proof}

\begin{proof}[Proof of Theorem \ref{thm:q_x-sinkhorn-descent}.]
Fix $k \in \NN$. It follows from \eqref{eq:alg_eq1} and Definition \ref{def:C_q-transform} that
\[
    \psi_{k-1}^{C_{q_\cdot}}(x) +  \int_\RRd \psi_{k-1}^*(a_k(x)+z) \,q_x(dz) =  a_k(x)\cdot x, \qquad \mu\text{-a.e. }x \in \supp(\mu).
\]
Then, by integrating both sides of the preceding equality with respect to $\mu$ and using \eqref{eq:alg_eq2}, we obtain
\begin{equation}
\label{eq:intermediate_equality_1}
    \int_\RRd \psi_{k-1}^{C_{q_\cdot}} \, d\mu + \int_\RRd \psi_{k-1}^* \, d\rho_k = \int_\RRd a_k(x)\cdot x \,\mu(dx).
\end{equation}
Similarly, \eqref{eq:alg_eq1} yields 
\[
    \psi_k^{C_{q_\cdot}}(x) = a_{k+1}(x)\cdot x - \int_\RRd \psi_k^*(a_{k+1}(x)+z) \,q_x(dz) \geq a_k(x)\cdot x - \int_\RRd \psi_k^*(a_k(x)+z) \,q_x(dz),
\]
for $\mu$-a.e. $x \in \supp(\mu)$, and integrating both sides of the preceding inequality with respect to $\mu$ gives
\begin{equation}
\label{eq:inequality1_dual_proof}
    \int_\RRd \psi_k^{C_{q_\cdot}} \, d\mu + \int_\RRd \psi_k^* \, d\rho_k \geq \int_\RRd a_k(x)\cdot x \,\mu(dx) = \int_\RRd \psi_{k-1}^{C_{q_\cdot}} \, d\mu + \int_\RRd \psi_{k-1}^* \, d\rho_k,
\end{equation}
where the last equality follows from \eqref{eq:intermediate_equality_1}.

Since $\psi_k$ is chosen so that $\psi_k$ attains the dual maximum covariance problem between $\rho_k$ and $\nu$, we have
\begin{equation}
\label{eq:dual_MCov}
    \MCov(\rho_k,\nu)=\inf_{\psi \text{ convex}} \int_\RRd \psi^* \, d\rho_k + \int_\RRd \psi \, d\nu=\int_\RRd \psi_k^* \, d\rho_k+\int_\RRd \psi_k \, d\nu,
\end{equation}
which implies that
\begin{equation}
\label{eq:inequality2_dual_proof}
    \int_\RRd \psi_k^* \, d\rho_k + \int_\RRd \psi_k \, d\nu \leq \int_\RRd \psi_{k-1}^* \, d\rho_k + \int_\RRd \psi_{k-1} \, d\nu.
\end{equation}
Thus,  we can conclude the desired monotonicity of the sequence $(\Ecal^{q_\cdot}(\psi_k))_{k \in \NN}$ by combining \eqref{eq:inequality1_dual_proof} and \eqref{eq:inequality2_dual_proof}
\begin{equation}
\label{eq:inequality3_dual_proof}
    \begin{aligned}
    \Ecal^{q_\cdot}(\psi_k) &= \int_\RRd \psi_k \, d\nu - \int_\RRd \psi_k^{C_{q_\cdot}} \, d\mu = \left( \int_\RRd \psi_k \, d\nu + \int_\RRd \psi_k^* \, d\rho_k \right) - \left( \int_\RRd \psi_k^{C_{q_\cdot}} \, d\mu + \int_\RRd \psi_k^* \, d\rho_k \right) \\
    &\leq \left( \int_\RRd \psi_{k-1} \, d\nu + \int_\RRd \psi_{k-1}^* \, d\rho_k \right) - \left( \int_\RRd \psi_{k-1}^{C_{q_\cdot}} \, d\mu + \int_\RRd \psi_{k-1}^* \, d\rho_k
    \right) = \Ecal^{q_\cdot}(\psi_{k-1}).
\end{aligned}
\end{equation}

Assume now that equality holds and $\rho_k\in\Prob_{p_2}(\RRd)$ and $\nu\in\Prob_{p_1}(\RRd)$ for some Hölder conjugates $p_1,p_2\in[1,\infty]$. Then equality holds in \eqref{eq:inequality2_dual_proof}, so that $\psi_{k-1}$ also attains the dual minimum in \eqref{eq:dual_MCov}. By Hölder's inequality, the maximum covariance problem between $\rho_k$ and $\nu$ has finite value. Moreover, the cost $c(z,y)=-z\cdot y$ satisfies the assumptions of \cite[Theorem 5.30]{Vi09}. Indeed, the required integrable lower bound follows from Young's inequality (or compact support if either $p_1 = \infty$ or $p_2 = \infty$). Hence the optimal coupling from $\rho_k$ to $\nu$ is unique and induced by the gradient of a convex function. Therefore
\[
    \nu=(\nabla\psi_{k-1}^*)_\#\rho_k.
\]
Moreover, the first-order optimality condition associated with the maximization in \eqref{eq:alg_eq1} gives
\[
    x = \int_\RRd \nabla\psi_{k-1}^*(a_k(x)+z) \,q_x(dz), \qquad \mu\text{-a.e. }x.
\]
Thus, if $M_0\sim\mu$, $Q|M_0=x\sim q_x$, and $M_1:=\nabla\psi_{k-1}^*(a_k(M_0)+Q)$, then $M_1\sim\nu$ and $M_0=\EE[M_1|M_0]$. Hence $\Law(M_0,M_1)$ is a $q_x$-Bass martingale from $\mu$ to $\nu$ generated by $v:=\psi_{k-1}^*$ and $a_k$. The optimality follows from Theorem \ref{thm:q_x-bass-optimality}.
\end{proof}

\begin{proof} [Proof of Proposition \ref{prop:q_x-sinkhorn-compactness}]
We first prove by induction that the sequence of potentials $(\psi_k)_{k \in \NN}$ is well-defined. Let $\psi=\psi_{k-1}$ and assume that $\psi \in L^1(\nu)$ is proper l.s.c. convex and that $\psi^*$ is Lipschitz. Define $G, H:\supp(\mu) \times \RRd \to \overline\RR$
\begin{equation}
\label{eq:G&H_def}
    G(x,a) := \int_\RRd \psi^*(a+z) \,q_x(dz), \qquad H(x,y):=(G(x,\cdot))^*(y)
\end{equation}
and denote by $L$ the Lipschitz constant of $\psi^*$. Then, for any $(x, a) \in \supp(\mu)\times\RRd$, \eqref{eq:G&H_def} implies
\[
    |G(x,a)| \leq |\psi^*(0)| + L|a| + L\int_\RRd|z| \,q_x(dz) < \infty,
\]
since $q_x \in \Prob_1(\RRd)$.
Moreover, for every $\eta\in\Pi(q_{x_1},q_{x_2})$ and $(x_1, a_1), (x_2, a_2) \in \supp(\mu)\times\RRd$,
\[
    |G(x_1,a_1)-G(x_2,a_2)|\leq L|a_1-a_2|+L\int_{\RRd\times\RRd}|z_1-z_2| \,\eta(dz_1,dz_2),
\]
and hence
\[
    |G(x_1,a_1)-G(x_2,a_2)|\leq L|a_1-a_2|+L\Wcal_1(q_{x_1},q_{x_2}).
\]
Therefore, the map $G$ is jointly continuous. Since $\mean(q_x)=x$, Jensen's inequality gives $G(x,a)\geq\psi^*(a+x)$ and therefore
\begin{equation}
\label{eq:upper_bound_H}
    H(x,y) = \sup_{a\in\RRd} \bigl(a\cdot y-G(x,a)\bigr) \leq \sup_{a\in\RRd}\bigl(a\cdot y-\psi^*(a+x)\bigr) = \psi(y) - x\cdot y.
\end{equation}
Since $\psi\in L^1(\nu)$ is convex, then $I := \interior(\co(\supp(\nu)))\subseteq\dom(\psi)$. Hence, \eqref{eq:upper_bound_H} implies that $H(x,\cdot)$ is finite on $I$, for any $x \in \supp(\mu)$. Furthermore, the inclusion $\supp(\mu) \subseteq I$ yields $\partial_y H(x,x) \neq \emptyset$ for every $x \in \supp(\mu)$. By conjugacy,
\[
    a \in \partial_y H(x,x) \quad \Longleftrightarrow \quad a \in \operatorname*{arg \,max}_{b\in\RRd} \bigl(b \cdot x - G(x,b) \bigr).
\]
We now record the following bound, which will also be used below. Choose $\delta>0$ such that
\[
    K_\delta := \{y \in \RRd: \operatorname{dist}(y,\supp(\mu)) \leq \delta\} \subseteq \interior(\co(\supp(\nu))).
\]
Since $\psi$ is locally bounded on $\interior(\co(\supp(\nu)))$, \eqref{eq:upper_bound_H} gives an upper bound for $H$ on $\supp(\mu)\times K_\delta$. Moreover,
\[
    H(x,y)\geq-G(x,0),\qquad |G(x,0)|\leq|\psi^*(0)|+L\int_\RRd|z| \,q_x(dz),
\]
and the $\Wcal_1$-continuity of $x\mapsto q_x$ together with the compactness of $\supp(\mu)$ gives
\[
    \sup_{x\in\supp(\mu)}\int_\RRd|z| \,q_x(dz)<\infty.
\]
Thus $H$ is bounded on $\supp(\mu)\times K_\delta$. If $a\in\partial_yH(x,x)$ and $|e|=1$, then
\begin{equation}
\label{eq:a_bounded_estimate}
    H(x,x)-H(x,x-\delta e)\leq\delta a\cdot e\leq H(x,x+\delta e)-H(x,x),
\end{equation}
so $|a|$ is bounded uniformly in $x$. Since $(x,a)\mapsto a\cdot x-G(x,a)$ is continuous, \cite[Theorem~18.19]{Ch06} yields a bounded Borel maximizer $a_k:\supp(\mu)\to\RRd$, and hence $(a_k)_\#\mu$ is compactly supported.

Define $\rho_k := (a_k(x)+z)_\#\bigl(q_x(dz)\mu(dx)\bigr)$. The boundedness of $a_k$ and the uniform first-moment bound on $(q_x)_{x\in\supp(\mu)}$ imply $\rho_k \in \Prob_1(\RRd)$. Since $q_x\ll\lambda^d$ for every $x\in\supp(\mu)$, also $\rho_k \ll \lambda^d$. Brenier's theorem therefore yields a proper l.s.c.\ convex function $\psi_k$ such that $\nu=(\nabla\psi_k^*)_\#\rho_k$.

Since $\nu$ is compactly supported, $\psi_k^*$ can be chosen globally Lipschitz. Moreover, by Fenchel--Young's equality,
\[
    \psi_k(\nabla\psi_k^*(x))=x\cdot\nabla\psi_k^*(x)-\psi_k^*(x) \in L^1(\rho_k),
\]
and therefore $\psi_k \in L^1(\nu)$. This closes the induction and proves that Algorithm~\ref{alg:q-sinkhorn} is well defined at every iteration.

Assume now that $(\mu,\nu)$ is irreducible, and fix a martingale coupling $\widehat\pi \in \Mart(\mu,\nu)$ satisfying \eqref{eq:irreducible_coupling}. Let $p_k \in \partial \psi_k(\bar x)$ and $\ell_k$ be as in the statement, so that $\widehat\psi_k\geq 0$ and $\widehat\psi_k(\bar x)=0$. Additionally,
\[
    \widehat\psi_k^*(a)=\psi_k^*(a+p_k)-p_k\cdot\bar x+\psi_k(\bar x), \qquad \widehat\psi_k^*(0)=0.
\]
Moreover, for any affine function $\ell$, we have $(\psi+\ell)^{C_{q_\cdot}}=\psi^{C_{q_\cdot}}+\ell$. Since $\mu$ and $\nu$ have the same mean, $\Ecal^{q_\cdot}$ is invariant under affine transformations. Therefore, the normalization introduced in the statement leaves the dual objective unchanged, and Theorem~\ref{thm:q_x-sinkhorn-descent} yields
\begin{equation}
\label{eq:strict_descent_normalized}
    \Ecal^{q_\cdot}(\widehat\psi_k)=\Ecal^{q_\cdot}(\psi_k)\leq\Ecal^{q_\cdot}(\psi_0).
\end{equation}

Since Jensen's inequality and $\mean(q_x)=x$ yield
\begin{align*}
    \widehat\psi_k^{C_{q_\cdot}}(x) & = \sup_{a\in\RRd} \left( a\cdot x -\int_\RRd \widehat\psi_k^*(a+z) \,q_x(dz) \right) \leq \sup_{a\in\RRd} \left( a \cdot x- \widehat\psi_k^*(a+x) \right)\\
    &= \sup_{b\in\RRd} \left( (b-x) \cdot x - \widehat\psi_k^*(b) \right) = \widehat\psi_k(x) - |x|^2,
\end{align*}
we have
\begin{equation}
    \label{eq:estimate_dual_obj}
    \begin{aligned}
        \Ecal^{q_\cdot}(\widehat \psi_k) & = \int_\RRd \widehat\psi_k(x) \, \nu(dx) - \int_\RRd \widehat\psi_k^{C_{q_\cdot}}(x) \, \mu(dx)\\
        & \geq \int_\RRd \widehat\psi_k(x) \, \nu(dx) - \int_\RRd \widehat\psi_k(x) \, \mu(dx) + \int_\RRd |x|^2 \mu(dx).
    \end{aligned}
\end{equation}

Therefore, combining \eqref{eq:strict_descent_normalized} and \eqref{eq:estimate_dual_obj} gives
\[
    \int_\RRd \widehat\psi_k(x) \, \nu(dx) - \int_\RRd \widehat\psi_k(x) \, \mu(dx) \leq\Ecal^{q_\cdot}(\psi_0)-\int_\RRd|x|^2 \,\mu(dx),
\]
which implies that 
\[
\sup_{k \in \NN} \int_\RRd \left( \int_\RRd \widehat \psi_k(y) \widehat \pi_x - \widehat \psi_k(x)\right) \mu(dx) < \infty.
\]
This means that we can invoke \cite[Lemma 3.9]{HaJoLoObPa26} to conclude that $(\widehat\psi_k)_{k\in\NN}$ is locally uniformly bounded on $I$.

Since $\nu$ is compactly supported, the Brenier potentials $\psi_k^*$ can be chosen with a common Lipschitz constant for all $k \in \NN$, which, for simplicity, we denote again by $L<\infty$. Therefore, since $\widehat\psi_k^*(0)=0$, this gives $|\widehat \psi_k^*(a)| \leq L|a|$ for every $a\in\RRd$ and $k \in \NN$. 

For any $k \in \NN$, set $\widehat a_{k+1} := a_{k+1} - p_k$ and define the functions $\widehat G_k, \widehat H_k:\supp(\mu)\times \RRd \to \RR$
\[
    \widehat G_k(x,a):=\int_\RRd \widehat \psi_k^*(a+z) \,q_x(dz),\qquad \widehat H_k(x,y):=(\widehat G_k(x,\cdot))^*(y).
\]
Then $\widehat a_{k+1}(x) \in \partial_y \widehat H_k(x,x)$. Now, note that on the same set $K_\delta$, we have the same estimates as above:
\[
    \widehat H_k(x,y) \leq \widehat \psi_k(y) - x\cdot y, \qquad \widehat H_k(x,y) \geq - \widehat G_k(x,0), \qquad |\widehat G_k(x,0)| \leq L \int_\RRd |z| \,q_x(dz).
\]
Hence, since the sequence $(\widehat \psi_k)_{k\in\NN}$ is locally uniformly bounded, $(\widehat H_k)_{k \in \NN}$ is uniformly bounded on $\supp(\mu)\times K_\delta$. Applying the estimate \eqref{eq:a_bounded_estimate} gives
\[
    \sup_{k \in \NN} \| \widehat a_{k+1}\|_{L^\infty(\mu)} < \infty.
\]

Finally, let $(k_n)_{n\in\NN}$ be an arbitrary subsequence. Since $(\widehat\psi_k^*)_{k\in\NN}$ is locally uniformly bounded and equi-Lipschitz, the Ascoli--Arzelà theorem yields, on every compact subset of $\RRd$, a uniformly convergent subsequence. Applying this on the sequence of compact sets $(\overline{B(0, k_n)})_{n \in \NN}$ and using a diagonal argument, we obtain a further subsequence, still denoted by $(k_n)_{n\in\NN}$, and a function $\varphi:\RRd\to\RR$ such that
\[
    \widehat\psi_{k_n}^*\longrightarrow \varphi
\]
locally uniformly on $\RRd$. Since convexity and the $L$-Lipschitz property are preserved under locally uniform convergence, $\varphi$ is a finite convex $L$-Lipschitz function.
Hence $\widehat\psi_{k_n}$ epi-converges to $\widehat\psi := \varphi^*$ by \cite[Theorem 11.34]{RoWe09}. The local uniform boundedness of $(\widehat\psi_k)_{k \in \NN}$ implies that $\widehat \psi$ is finite on $I$, and therefore 
\[
    \widehat \psi_{k_n} \longrightarrow \widehat \psi \quad\text{locally uniformly on } I.
\]
Since $\varphi$ is finite and continuous, $\widehat\psi^*=\varphi$, which proves the claim.
\end{proof}

\printbibliography

\end{document}